\documentclass[preprint]{elsarticle}
\usepackage{array}
\usepackage{amsmath,amssymb}
\usepackage{breqn}
\usepackage{subcaption}
\usepackage{stackengine}
\usepackage{verbatim} 
\usepackage{algorithm}
\usepackage{algorithmic}
\newcolumntype{R}{>{\centering\arraybackslash}m{8cm}}
\usepackage{graphicx}
\usepackage{xcolor}
\usepackage[inline]{enumitem}
\usepackage{nomencl}
\makenomenclature

\usepackage{hyperref}

\DeclareMathOperator*{\Min}{Minimize}

\usepackage{etoolbox}
\renewcommand\nomgroup[1]{%
  \item[\bfseries
  \ifstrequal{#1}{A}{Sets}{%
  \ifstrequal{#1}{B}{Parameters}{%
  \ifstrequal{#1}{C}{Decision Variables}{%
  \ifstrequal{#1}{D}{Functions}{}}}}%
]}

\newtheorem{theorem}{Theorem}

\newproof{proof}{Proof}

\begin{document}

\begin{frontmatter}

\title{Modeling COVID-19 optimal testing strategies in long-term care facilities: An optimization-based approach }

\author[a,b]{Mansoor Davoodi\corref{cor1}}
\ead{m.davoodi-monfared@hzdr.de}

\author[a,b]{Ana Batista}
\ead{a.batista-german@hzdr.de}

\author[a,b]{Abhishek Senapati}
\ead{a.senapati@hzdr.de}

\author[a,b]{Weronika Schlechte-Welnicz}
\ead{w.schlechte-welnicz@hzdr.de}

\author[c]{Birgit Wagner}
\ead{birgit.wagner@dwlz.de}

\author[a,b,d,e]{Justin M. Calabrese}
\ead{j.calabrese@hzdr.de}


\cortext[cor1]{Corresponding author.}
\address[a]{Center for Advanced Systems Understanding (CASUS), D-02826 G$\ddot{o}$rlitz, Germany}
\address[b]{Helmholtz-Zentrum Dresden-Rossendorf (HZDR), D-01328 Dresden, Germany}
\address[c]{Diakonisches Werk im Kirchenbezirk Löbau-Zittau GmbH}
\address[d]{Helmholtz Centre for Environmental Research-UFZ, Leipzig, Germany}
\address[e]{Dept. of Biology, University of Maryland, College Park, MD, USA}


\begin{abstract}
Long-term care facilities have been widely affected by the COVID-19 pandemic. Retirement homes are particularly vulnerable due to the higher mortality risk of infected elderly individuals. Once an outbreak occurs, suppressing the spread of the virus in retirement homes is challenging because the residents are in contact with each other and isolation measures cannot be widely enforced. Regular testing strategies, on the other hand, have been shown to effectively prevent outbreaks in retirement homes. However, high frequency testing may consume substantial staff working time, which results in a trade-off between the time invested in testing, and the time spent providing essential care to residents.  
Thus, developing an optimal testing strategy is crucial to proactively detect infections while guaranteeing efficient use of limited staff time in these facilities.  
Although numerous efforts have been made to prevent the virus from spreading in long-term care facilities, this is the first study to develop testing strategies based on formal optimization methods.
This paper proposes two novel optimization models for testing schedules. The models aim to minimize the risk of infection in retirement homes, considering the trade-off between the probability of infection and staff workload. We employ a probabilistic approach in conjunction with the optimization models, to compute the risk of infection, including contact rates, incidence status, and the probability of infection of the residents.
To solve the models, we propose an enhanced local search algorithm by leveraging the \textit{symmetry property} of the optimal solution. We perform several experiments with realistically sized instances and show that the proposed approach can derive optimal testing strategies. 
\end{abstract}


\end{frontmatter}

\section{Introduction}

Long-Term Care Facilities (LTCF) include institutions such as Retirement Homes (RH), nursing homes, and rehabilitation centers that provide medical and personal support to patients \cite{danis2020high}. These organizations aim to provide high-quality care for the elderly population and efficient management of resources \cite{anderson2021optimization, spasova2018challenges}. 

The coronavirus disease 2019 (COVID-19) pandemic has had a significant impact on LTCFs \cite{gmehlin2020covid, giri2021nursing}. These facilities typically have a high density of elderly people at a higher risk for mortality after being infected with the severe acute respiratory syndrome coronavirus 2 (SARS-CoV-2) virus.
%
The majority of residents in these institutions have pre-existing conditions (i.e., diabetes, respiratory disease, hypertension, chronic heart diseases) \cite{childs2019burden, moyo2020risk}, which have been linked to an increased risk of death in older patients \cite{birgand2021testing, roxby2020outbreak}. 
Early in the development of the COVID-19 pandemic, it became clear that the elderly are the most impacted. For instance, in the European Union/European Economic Area (EU/EEA), people over 65 years old accounted for 88\% of all COVID-19 related deaths. In particular, LTCFs have been linked to 37-62\% of fatalities in several EU/EEA countries \cite{comas2020mortality, european2020surveillance}. Similarly, in the United States, over 30\% of COVID-19 related deaths were associated with nursing homes \cite{centers2020covid}.
%

Although COVID-19 vaccines have demonstrated high efficacy, they do not provide full immunity against infection nor do they completely mitigate mortality risk in older individuals  \cite{thomas2021safety, nasreen2021effectiveness, kahn2021mathematical, patriarca1985efficacy}.
Thus, to contain the rapid spread of the virus in LTCFs, the Centers for Disease Control (CDC) and the European Centre for Disease Prevention and Control (ECDC) have issued infection prevention and control (IPC) recommendations, including social distancing, daily screening (testing) of staff and residents, isolation, and visitation restrictions \cite{gmehlin2020covid, danis2020high}. Among the IPC recommendations, widespread testing of staff and residents has been demonstrated to detect infections proactively and avoid propagation.
For instance, \cite{roxby2020outbreak} performed an outbreak investigation in a nursing home in Seattle and found that symptom screening fails to identify infected residents. In contrast, preventive testing combined with safety strategies can reduce viral spread. Similarly, \cite{telford2020preventing}, and \cite{graham2020sars} demonstrated that facility-wide testing is an efficient strategy to prevent COVID-19 outbreaks since it helps to identify asymptomatic infections in LTCFs proactively. Cohort isolation of positive residents in conjunction with widespread screening has also been considered as an efficient strategy to prevent the spread of the virus in \cite{krone2021control}.
Overall, widespread testing in LTCFs has been shown to be one of the most efficient strategies to limit outbreaks. The proactive testing strategy combined with other protocols, such as isolation and quarantine, facilitates the timely implementation of control procedures due to the early detection of the virus \cite{kucharski2020effectiveness, bergstrom2020frequency}.

Despite recommendations from health authorities to mitigate the spread of COVID-19, LTCFs are still at risk of outbreaks. Besides the factors related to the residents' health, other aspects have been reported as the main causes of outbreaks, including understaffing, residents sharing common spaces, and high contact rates between residents and staff \cite{gorges2020staffing, giri2021nursing, chen2021nursing}. In particular, RH facilities have faced difficulties in controlling the spread of the virus because isolating infected residents is often not possible due space limitations and regulations. Moreover, there is a shortage of staff to administer the test to the residents \cite{onlineRHreport, krone2021control}. 
As a consequence, once an infection arrives at the facility, it spreads very rapidly, putting both residents and staff at risk. In order to prevent the pandemic from spreading to vulnerable populations, it is therefore necessary to implement effective testing strategies in RHs.

Performing testing procedures in RHs is a challenge. 
Firstly, the staff, who have a defined workload of caregiving tasks, need to dedicate a portion of their work time to administering tests to residents, which may affect the quality of care services. Due to financial and/or regulatory limitations, hiring additional staff during the pandemic to perform testing activities might not be possible \cite{onlineRHreport}. Secondly, implementing frequent testing may cause discontent among the residents because of the implications of the test procedure. These challenges evidence a fundamental trade-off between the staff workload for care duties and testing the residents.
%
An \textit{optimal testing strategy} therefore needs to satisfy staff workload limitations while minimizing the risk of infection in the facilities.

%

Motivated by the current challenges in preventing the spread of COVID-19 in RHs, we propose a novel optimal testing strategy for outbreak suppression, in the framework of an optimization model. The testing procedure in RHs is usually performed by trained staff. The residents are divided into groups and tested considering a predefined test schedule. The staff is responsible for cleaning, disinfecting, preparing the testing workspace, and administering tests, all of which consume a significant amount of time. 
%
Therefore, finding an optimal testing schedule requires determining the testing interval, the number and size of groups, and the day on which to test each group. Collectively, these considerations result in a challenging combinatorial optimization problem. In the Operations Research literature, this problem is similar to the resource allocation problem that optimally assigns resources to activities to minimize related costs. Due to the complexity of solving these combinatorial problems, heuristic solutions are mainly considered, such as search algorithms. These algorithms have guaranteed efficient performance for obtaining global solutions \cite{rios2021location, bouajaja2017survey, moreira2015model}. 
To the best of our knowledge, there are no studies in the literature that focus on developing testing strategies in LTCFs that take formal optimization methods into account.

The contribution of this paper is twofold. First, we introduce novel optimization models for testing schedule strategies. Specifically, we develop two Mixed Integer Nonlinear Programming (MINLP) models for balancing the staff's workload in RHs while minimizing the expected detection time of a probable infection inside the facility. The first model minimizes the expected time to detect an infection, considering a threshold on the maximum portion of staff time allocated for the testing process. The second model minimizes the testing workload for the staff, incorporating the number of infections in the neighborhood. In both cases, the expected risk of infection is computed via a probabilistic disease transmission model.
Second, to solve the models, we propose a highly useful property --which we denote the \textit{Symmetry property}-- and leverage it to propose an enhanced local search algorithm able to find optimal solutions.


This paper is organized into five sections. After reviewing related and recent studies in Section \ref{RelatedWork}, the models for the problem of finding optimal testing strategies in RHs are developed in Section \ref{Modeling}. A practical approach to finding the optimal testing strategy is proposed in Section \ref{SolutionApproach}. Simulation results for different scenarios are presented in Section \ref{Simulation}. Finally, a conclusion is drawn in Section \ref{Conclusions}.
\section{Related Work}
\label{RelatedWork}

This section covers the literature related to the optimization of testing strategies during the COVID-19 pandemic. We focus on works that study the design of testing strategies for LTCFs aiming to control the spread of the virus. 
 

Developing an effective testing strategy is crucial to prevent the spread of the virus that causes COVID-19. However, a limited number of studies have considered models within a mathematical optimization framework aiming to provide optimal testing strategies to control infectious diseases in LTCFs (See \cite{jordan2021optimization, choi2021fighting} for a review of Optimization in the context of COVID-19). Most studies have focused on studying testing strategies employing compartmental models (i.e., Susceptible- Infectious-Recovered (SIR) and variants) \cite{kermack1927contribution} and simulation models \cite{jordan2021optimization}.


Viral testing has different aims depending on the application context \cite{yuval2020optimizing}: \begin{enumerate*}
    \item  Diagnostic testing: testing of symptomatic patients or those who had contact with infected individuals.  
    \item  Spread suppression: widespread testing of asymptomatic individuals. 
    \item  Outbreak detection: randomized testing of asymptomatic individuals when disease prevalence is low. 
\end{enumerate*}
Several studies have been dedicated to the optimization of diagnostic testing strategies (See e.g., \cite{dhiman2021adopt, elaziz2020new}). However, the main focus of diagnostic testing is clinical, for detecting and treating the presence of the virus in individuals \cite{prinzi2020screening}, and it is out of the scope of our study. We review studies addressing testing strategies for outbreak detection and spreading suppression.

Some studies have considered outbreak detection strategies in nursing homes. The authors in \cite{yuval2020optimizing} remarked that outbreak detection strategies are recommended for small institutions with a limited budget of testing and in which new infections are rare (i.e., low transmission rates). They developed a compartmental network-based Susceptible-Exposed- Infectious-Recovered (SEIR) model considering the heterogeneity of connections, incubation period, and test efficacy. They found that testing small groups with high frequency is a better strategy for outbreak detection than testing larger groups less frequently. 
%
In Germany, \cite{german2021bulletin} studied COVID-19 outbreaks in retirement homes and found that symptom control and testing, in addition to vaccination, are effective prevention strategies. \cite{lucia2021modeling} employed an SEIR network epidemic model based on disease status. They used the model to study a shield-immunity approach, considering a bipartite network between the staff and residents. The results showed that shield-immunity interventions in conjunction with regular testing helps to reduce the size of the outbreak. 
\cite{see2021modeling} considered an outbreak testing strategy employing a Reed-Frost model. The authors evaluated the effectiveness of outbreak testing for staff and residents in nursing homes. Their findings suggested that combining infection control practices with massive testing is an effective approach to prevent the spread of COVID-19. 

In the context of spreading suppression strategies in LTCFs, \cite{tsoungui2021preventing} studied testing strategies for closed facilities (i.e., LTCFs and incarceration centers). The study considers an extended deterministic SEIR model to evaluate the impact of widespread testing on the staff on the number of resident infections. The results showed a 40\% reduction in cases by minimizing the number of contacts between staff and residents and testing the staff every five days. The authors remarked that these results are highly dependent on the type of facility. 
The authors in \cite{bagger2022reducing} studied strategies for reopening activities considering a Markov process. They considered a graph representation of the individuals' contacts to determine the structures that reduced the spread of the disease. To reduce disease spread, their results suggeste that limiting the interactions of participants in an activity is more effective than a size limitation.

Within the general COVID-19 context, most studies have focused on developing personnel scheduling models for healthcare workers \cite{guler2020decision, seccia2020nurse, zucchi2021personnel, sanchez2021modelling, jiang2021optimal, guerriero2021modeling}, but disregarding an epidemiological model for computing the risk of infection. Of particular interest is the work of \cite{abdin2021optimization}, where the authors developed an optimization-based compartmental model for planning, testing, and control. Similarly, \cite{calabrese2022optimal} employed an SEIR model to study the optimal balance between spreading suppression and outbreak detection testing strategies under limited testing capacities.
However, the solution framework of these studies is focused on country-wise strategies and not on particular organizations. Overall, to the best of our knowledge, the problem of developing a testing schedule that formally optimizes the trade-off between staff workload and its impact on the number of infections in LTCFs has not been studied in the literature. 


\section{Modeling the test scheduling problem for retirement homes}
\label{Modeling}

The required notation for the proposed models is as follows:
 
\nomenclature[B01]{$m$}{Number of residents.}
\nomenclature[B02]{$n$}{Number of staff.}
\nomenclature[B03]{$P_{time}$}{Preparation time for cleaning and preparing the testing workspace for each group of residents.}
\nomenclature[B04]{$T_{time}$}{Testing time for each resident.}
\nomenclature[B05]{$p$}{ Maximum portion of staff time which can be allocated to testing processes\footnote{Note that $p$ is a decision variable in Model 2.}}
\nomenclature[B06]{$\alpha$}{Acceptable level of risk inside the retirement home compared to the outside (background~risk).}
\nomenclature[B09]{$Max_{\tau}$}{Upper bound on the test interval.}
\nomenclature[B10]{$\beta$}{Probability of disease transmission per one contact.}
\nomenclature[B11]{$\kappa$}{Average number of daily contacts for each resident.}

\nomenclature[C1]{$k$}{ Number of groups of residents.}
\nomenclature[C2]{$\tau$}{Testing interval.}
\nomenclature[C3]{$G=\{g_1,g_2,\dots,g_k\}$}{Size of groups, i.e., $g_i$ denotes the size of $i^{th}$ group $k$.}
\nomenclature[C4]{$D=\{d_1,d_2,\dots,d_k\}$}{Day for testing each group $k$.}

\printnomenclature

In the problem of scheduling the residents for testing in RHs, the aim is to find an optimal test schedule that minimizes the risk of infection in the facility while balancing the staff workload. Thus, the testing strategy includes defining a testing interval, the testing day, and the grouping of residents. Also, the tests are available and there is no budget constraint to obtain them\footnote{The problem discussed in this paper and computational experiments are based on data of the retirement homes run by Diakonisches Werk im Kirchenbezirk Löbau-Zittau GmbH, Saxony, Germany.}.

The RH consists of $m$ residents and $n$ staff. For simplicity, we assume homogeneity in terms of the probabilities of transmitting the virus and being infected. Indeed, in the model, we utilize the average values of such parameters. We also assume that, due to strict governmental regulations, the staff is subject to regular testing, so we only focus on testing the residents in the facility \cite{european2020surveillance}. 
%
%

The testing process is described as follows. The residents are divided into $k$ groups, and tested in predefined time intervals, such that every $\tau$ days one test per resident is performed. Thus, a testing strategy is defined using a quadruple $(k,\tau, G,D)$ where $k$ is the number of groups for testing, $\tau$ is the testing interval, $G=\{ g_1,g_2,\dots,g_k\}$ is a partitioning of the residents to $k$ groups, and $D=\{d_1,d_2,\dots,d_k\}$ shows the testing day for each group. The $m$ residents are divided into $k$ groups, and the test is performed in a day $d_i$, where $0 < d_i \leq \tau$, for $i=1,2,\dots,k$. Note that all decision variables $k$ $\tau$, $G$, and $D$ are defined in the integer domain.
Without loss of generality, we set the reference day as zero, so $d_k=\tau$. Each group is tested in one batch for the staff, who cleans and prepares the testing workspace. We denote the preparation cost (here time or workload) by $P_{time}$. In addition to this cost, each resident has his/her own testing time, denoted by $T_{time}$. For simplicity, we consider the costs in terms of one working day. The total cost for one round of testing of all residents is expressed as follows:

\begin{equation}
Testing~cost = k \times P_{time} + m \times T_{time}.
\label{Testing_Cost}
\end{equation}

The testing cost in Equation \eqref{Testing_Cost} is the total time the staff spend in preparing and performing tests, every $\tau$ days. Model 1 is a MINLP defined as follows: 

\textbf{Model 1:} 
\begin{align} \label{M1.2}
& \Min ~~~ \text{Expected Detection Time of} ~~ (k,\tau, G,D)\\ \label{M1.3}
& \text{s.t.:}   \notag \\ 
& k \times P_{time} + m \times T_{time} \leq p \times n \times \tau\\ \label{M1.4}
&\sum_{i=1}^{k} |g_i| = m\\  \label{M1.5}
& \tau \leq Max_{\tau}\\  \label{M1.6}
&|g_i| \leq Max_g,~~~~ \forall i=1,2,\dots, k 
\end{align}


The objective function \eqref{M1.2} aims to minimize the expected detection time defined by the quadruple $(k,\tau, G,D)$. For any testing strategy $(k,\tau, G, D)$, it is possible to compute the expected time to detect a probable infection among the residents. For example, if $k=1$, the expected time to detect will be $\frac{\tau}{2}$, and for $k>1$, it depends on the rate of COVID-19 transmission and the number of contacts among the residents. At the end of this section, we discuss how to compute the expected detection time (See Subsection \ref{Computingdetection}).

Constraint \eqref{M1.3} limits the proportion of staff time dedicated to the testing process, where $p$ refers to the portion of the staff's working time which can be allocated to testing. Constraint \eqref{M1.4}, imposes that the sum of the defined groups is equal to the total number of residents, $m$, where $|g_i|$ denotes the size of the group $g_i$.   
Constraint \eqref{M1.5} ensures that the testing interval, $\tau$, does not exceed a predefined upper bound, $Max_{\tau}$. For instance, $Max_{\tau}=7$ forces the algorithm to find a solution that guarantees each resident is tested at least once a week. 
Constraint \eqref{M1.6}, limits the size of each batch to the upper bound,  $Max_g$, of the maximum number of groups. 

The resulting strategy of Model 1 aims to find the testing strategy that minimizes the time for detecting a probable infection in the facility under the set of constraints \eqref{M1.3}--\eqref{M1.6}. We show under the homogeneity assumption of residents, it also results in minimizing the expected number of infections. 

\begin{theorem}
In the problem of finding an optimal testing strategy for a retirement home assuming homogeneity of residents, the objective of minimizing the detection time is equivalent to minimizing the number of infections after an infection arrives at the retirement home.
\label{thm_eq}
\end{theorem}

See proof \ref{thm_eq_proof} in the Appendix.

In Model 1 (i.e., Equations \eqref{M1.2}--\eqref{M1.6}), we derive an optimal testing strategy by minimizing the expected detection time. Alternatively, a manager of an RH might also want to reduce the risk of infection in the facility to an admissible preferred level. Then he/she wishes to achieve this level with the minimum possible allocation of staff time. In order to represent this problem, we define Model 2 as follows:

\textbf{Model 2:} 
\begin{align} \label{M2.7}
& \Min ~~~ p \\ \label{M2.8}
& \text{s.t.:}   \notag \\ 
&Pr(infection|(k,\tau, G,D)) \leq \alpha \times background~risk\\ \label{M2.9}
&k \times P_{time} + m \times T_{time} \leq p \times n \times \tau\\ \label{M2.10}
&\sum_{i=1}^{k} |g_i| =m\\  \label{M2.11}
& \tau \leq Max_{\tau}\\  \label{M2.12}
& |g_i| \leq Max_g,~~~~ \forall i=1,2,\dots, k.
\end{align}

The objective function \eqref{M2.7} minimizes $p$, the portion of the staff time dedicated to the testing process. Constraints \eqref{M2.9} -- \eqref{M2.12} are as defined in Model 1. 
Constraint \eqref{M2.8}, guarantees that the probability of infection per resident under a testing strategy $(k,\tau, G, D)$ does not exceed a desired level of infection, $\alpha$. So that $\alpha>0$ is a coefficient defined by the manager to set the level of risk inside the RH.
In order to measure the magnitude of such a level, we introduce the \textit{background risk}, that is, the probability of infection for an individual who lives in the area of the focal RH. It is a reference point to help the manager compare the risk of infection inside the RH with that in the surrounding region. 
We remark that the background risk has been defined to determine the probability that an infection arrives at the facility either via staff or visitors. The background risk can be computed using the number of incidences in the local neighborhood.  

In the following subsections, we describe the probabilistic approach developed to compute the expected detection time of an infection. 

\subsection{Computing the probability of infection}
\label{Computingprobability}

This study aims to minimize the probability of an outbreak spreading among the residents. This section presents a probabilistic approach to compute the risk of infection in a RH. 

We assume that a resident can be infected at time $t=0$, so we compute the probability of infection for a resident at any time $t=d$. In order to derive the probability of infection, we define the following parameters:

\begin{itemize}
    \item $\beta$: The probability of transmission from one infected resident to a susceptible resident per contact.
    \item $\kappa$: The average number of contacts per resident. 
\end{itemize}

Let $P_I(m,\kappa,\beta,d)$ be the probability of infection for a resident at day $t=d$ in a RH with $m$ residents who have $\kappa$ contacts with other residents per day. Let $s$ be the source (first resident) of infection at day $t=0$, and $u$ be a resident who stays healthy till day $t=d-1$. There are two ways $u$ gets infected on day $d$; via some direct contact with $s$, or via contact with one of the other $m-2$ residents.
The transmission probability for each of the contacts is $\beta$, and the probability of infection for the source is one, while the probability of infection for the other residents is $P_I(m,\kappa,\beta,d-1)$. That means that the probability of an arbitrary resident like $u$ staying healthy after $c$ contacts with $s$ is $(1-1\times \beta)^c$. While the probability of $u$ staying healthy after $c$ contacts with the other $m-2$ resident in day $d$ is $(1-P_I(m,\kappa,\beta,d-1) \times \beta)^c$. 
Therefore, the infection probability can be represented as a recursive equation. Given that we consider $\kappa$ contacts for each resident per day and all of them have an equal chance to occur, the  equation can be defined as follows:

\begin{dmath}
P_I(m,\kappa,\beta,d)=1-(1-P_I(m,\kappa,\beta,d-1)) \times \left[(1-\beta)^{\kappa \frac{1}{m-1}} \times (1-P_I(m,\kappa,\beta,d-1)\beta)^{\kappa (1-\frac{1}{m-1})} \right], 
\label{Prob_Infection}
\end{dmath}

where $\kappa \frac{1}{m-1}$ is the expected number of contacts between $u$ and $s$, and $\kappa (1-\frac{1}{m-1})$ is the expected number of contacts between $u$ and the other residents except $s$. Note that, $(1-P_I(m,\kappa,\beta,d-1))$ is the probability of $u$ staying healthy until day $t=d-1$.

Equation \ref{Prob_Infection} provides a simple recursive formula for the probability of infection per resident $d$ days after an infection first arrives at the facility. We can also easily take into account the role of staff as an intermediate node between two residents for indirect contacts (i.e., adding the number of contacts between residents and staff to $\kappa$). However, for simplicity, we keep the formula for the direct contact network among the residents and utilize it for computing the expected detection time for a testing strategy $(k,\tau, G, D)$.

\subsection{Computing the Expected detection time}
\label{Computingdetection}

\begin{figure}[H]
\centering
\includegraphics[width=12cm]{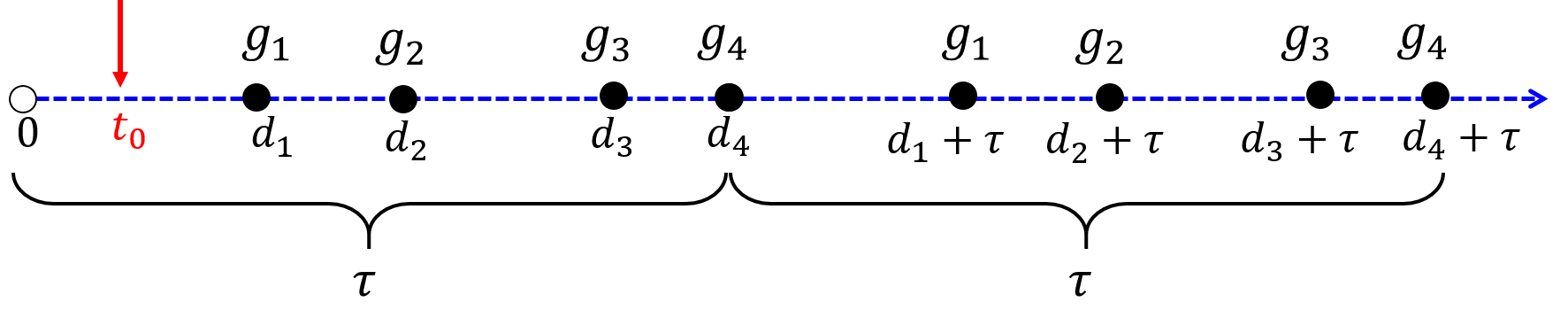}
\caption{Schematic view of a testing strategy with $k=4$ groups and two time periods of $\tau$.}
\label{Fig_strategy}
\end{figure}

In this subsection, we describe how to compute the expected detection time for a given testing strategy $(k,\tau, G, D)$.
Note that $G=\{g_1, g_2,\dots,g_k\}$ are the $k$ groups of residents and $D=\{d_1, d_2,\dots,d_k\}$ are the days on which they are tested. For the sake of simplicity, we assume an origin $t=0$ and define $d_i$ as the distance (i.e., number of days) from the origin, and $d_i \leq d_{i+1}$. Figure \ref{Fig_strategy} shows a configuration of such strategy with $k=4$ groups for two periods $\tau$.

Suppose an infection arrives at the RH at time $t=t_0$. Let $t_0 \leq d_1$, otherwise, we can reorder the groups because the configuration is circular. Since the residents are partitioned into $k$ groups, clearly the infection should be detected in one of the days $t=d_1,t=d_2,\dots ,t=d_k$. Thus, it is sufficient to compute the probability of detecting the infection on these predetermined days. To this end, we apply the probability function $P_I(m,\kappa,\beta,d)$ defined in Equation \ref{Prob_Infection} with some manipulations. The function $P_I(m,\kappa,\beta,d)$ is the probability of one resident getting infected $d$ days after the origin. So, if a group of $r$ randomly selected residents is tested after $d$ days, the probability of at least one of them being infected is 

\begin{equation}
Pr(m,\kappa,\beta,d,r)=\frac{r}{m}+ (1-\frac{r}{m})\times \left( 1-[1-P_I(m,\kappa,\beta,d)]^r \right),
\end{equation}

where $\frac{r}{m}$ is the probability of the source of infection (the first resident who gets infection) being sampled in the group. Let us define the complement of this probability function by $\bar{Pr}(m,\kappa,\beta,d,r)=1-Pr(m,\kappa,\beta,d,r)$, that is, the probability that the introduced infection is not detected after $d$ days by testing a group of size $r$. Now, we can compute the expected detection time for an introduced infection at day $t=t_0$ as

\begin{dmath}
\mathbb{E} (k,\tau, G,D,t_0)= \sum_{i=1}^{k} (d_i-t_0) \times \left [\left[ x +(1-x) \times Pr(m,\kappa,\beta,d_i-t_0,|g_i|) \right] \times \prod_{j=1}^{i-1} \bar{Pr}(m,\kappa,\beta,d_j-t_0,|g_j|) \right],
\label{Exp_detect_t0}
\end{dmath}

where $x=\frac{|g_i|}{m-\sum_{j=1}^{i-1} |g_j|}$ is the probability that the source of infection is sampled in $g_i$, and $\sum_{j=1}^{0} |g_j|$ is $0$. In this calculation, $\prod_{j=1}^{i-1} \bar{Pr}(m,\kappa,\beta,d_j-t_0,|g_j|)$ is the probability that the infection is not detected before day $d_i$. So, for convenience, we set $\prod_{j=1}^{0} \bar{Pr}(m,\kappa,\beta,d_j-t_0,|g_j|)=1$. Finally, the expected time for detecting an infection under a testing strategy $(k,\tau, G,D)$ is computed as follows:

\begin{equation}
Expected~Detection~ Time~ of~(k,\tau, G,D)= \frac{1}{\tau} \sum_{t_0=0}^{\tau-1} \mathbb{E} (k,\tau, G,D,t_0).
\label{Exp_time}
\end{equation}

Thus, for any testing strategy $(k,\tau, G,D,t_0)$ in a RH with $m$ residents, the expected detection time can be computed using Equation \ref{Exp_time}, that is the objective function of Model 1. If $\bar{\tau}$ is such an expected detection time, the probability of infection per resident in the day the infection is detected can be computed by Eq. \ref{Prob_Infection} for $d=\bar{\tau}$. That is, the value of $Pr(infection|(k,\tau,G,D))$ in Constraint \eqref{M2.10} of Model 2.

\section{An efficient approach for computing the optimal testing strategy}
\label{SolutionApproach}

Computing a testing strategy $(k,\tau,G,D)$ means determining the number of groups, $k$, a testing interval, $\tau$, the size of groups, $G=\{g_1,g_2,\dots,g_k\}$, and the testing day, $D=\{d_1,d_2,\dots,d_k\}$. The values are integer decision variables, and a feasible solution should be selected from the integer space.
Clearly, if the strategy is chosen from the real space, $\mathbb{R}$, the optimal solution will result in a better objective value, which is a lower bound on the objective value of the optimal solution when it is chosen from the integer space. 

This section presents the solution approach for solving the testing scheduling problem for RHs. We first present a useful observation --which we call the \textit{symmetry property}-- of the optimal solution. Then, we develop an enhanced local search algorithm for solving any instance of the problem based on the symmetry property.

\paragraph{The symmetry property} 
The symmetry property of the optimal solution indicates that for a given test interval $\tau$ and the number of groups $k$, the optimal solution can be obtained by evenly distributing the $m$ residents in the interval $\tau$. More precisely, $|g_1|=|g_2|=\dots=|g_k|=\frac{m}{k}$, and $d_i=i\frac{\tau}{k}$, for $i=1,2,\dots,k$. This property is similar to Purkiss's principle \cite{waterhouse1983symmetric} of the symmetry functions, but not the same because the objective function is not invariant under the possible permutations of pairs $(k,\tau)$.

We show the correctness of the symmetry property by providing mathematical proof for the extreme cases of the transmission rate, $\beta$, and the number of contacts, $\kappa$. Further, for the general case, we provide a brute force analysis on real-world sized instances of the problem. If a symmetry strategy is an integer solution, it will be the optimal strategy of the problem as well. Otherwise, we use the optimal solution as the seed of a global heuristic search to find the optimal integer strategy. Note that, for a given $k$ and $\tau$, the number of possible different combinations of $G$ and $D$ are of the order $O\left( \left( \begin{array}{c} \tau \\ k \end{array} \right)  \left( \begin{array}{c} m \\ k \end{array} \right) \right)$.

\begin{theorem}
The optimal testing strategy of the real space is a symmetry strategy for the cases $\beta \longrightarrow 0$, $\beta \longrightarrow 1$, $\kappa \longrightarrow 0$, or $\kappa \longrightarrow +\infty$.
\label{thm_limit}
\end{theorem}

The proof \ref{thm_limit_proof} is contained in the Appendix.

In the problem of finding the optimal testing strategy, the cases $\beta \longrightarrow 0$ and $\kappa \longrightarrow 0$ indicates very low risk of the disease propagating through the retirement home, and the cases $\beta \longrightarrow 1$ and $\kappa \longrightarrow +\infty$ indicate the opposite. 
Theorem \ref{thm_limit} shows that for both situations, the symmetry strategy is the optimal solution. We strongly believe this holds for any value of $\beta \in [0,1]$ and any contact number $\kappa \geq 0$. We investigated the correctness of this claim by sensitivity analysis for almost all real-world scales of the problem, that is, for any combination of parameters $k \leq \tau <15$ and $m \leq 100$. That means having at least one test per two weeks and a retirement home with fewer than 100 residents. We simulated the continuous space with the precision $0.01$ and, using a brute force algorithm, investigated the optimality of the symmetry strategy. However, due to the complexity of the function to compute the expected detection time, it is an intractable problem to show this fact. So, we left the mathematical proof as an open problem.

As previously mentioned, a feasible strategy for the RH testing schedule problem is a solution with the integer values $k$, $\tau$, and the integer sets $G$ and $D$. Hence, whether the symmetry property always holds or not, it is applicable only when integer values generate the obtained symmetry solution. Otherwise, with a high probability, the integer neighbors of the symmetry solution will be the optimal practical solution of the problem. Based on this fact, we propose a global heuristic local search algorithm. We remark that for a small size of the problem (e.g., $m<50$, $Max_g<20$, $n<10$ and $k \leq \tau \leq 7$) the computation for searching on all (integer) possible solutions of the problem and returning the optimal one is not intractable and can be implemented in practice. We used this approach to test the computational performance of the proposed algorithm in finding the optimal solution.

\paragraph{\textbf{Enhanced local search algorithm}}

The proposed enhanced local search algorithm, first generates all pairs of $(k,\tau)$ which satisfy Constraint \eqref{M1.2} of Model 1, $k \times P_{time} + m \times T_{time} \leq p \times n \times \tau$. For Model 2, since $p$ is a decision variable, the pairs of $(k,\tau)$ can be generated one-by-one, and the solutions with $p>1$ are infeasible. Then, the algorithm calculates the symmetry strategy, say $S$ for any feasible pair of $(k,\tau)$. If $S$ is an integer solution and satisfies the other constraints of the problem (e.g., $|g_i| \leq Max_g,~for~i=1,2,\dots, k$, or $Pr(infection|(k,\tau, G,D)) \leq \alpha \times background~risk$ in Model 2), the algorithm reports it as the optimal strategy, otherwise, a heuristic search on the integer possible solutions around $S$ is implemented.

The pseudocode of the algorithm for solving Model 1 is presented in Algorithm \ref{AlgorithmFirstModel} in the Appendix. The proposed heuristic search algorithm is similar to the \textit{Simulated Annealing} technique \cite{kirkpatrick1983optimization} but only on the integer neighborhoods of groups' size, $G$ (See Algorithm \ref{HeuristicSearch} in the Appendix). In order to efficiently implement this heuristic search, the algorithm utilizes the symmetry property on the feasible pairs of $(k,\tau)$. That is, for two feasible pairs $(k_1,\tau_1)$ and $(k_2,\tau_2)$, if $S_1=(k_1,\tau_1,G_1,D_1)$ is an integer feasible solution for the problem and its objective value is better (less detection time in Model 1, or smaller allocated workload for testing process in Model 2) than the symmetry strategy of $(k_2,\tau_2)$, then it prunes the case $(k_2,\tau_2)$ without searching on its possible integer neighbors. In fact, $(k_2,\tau_2)$ is not promising anymore, and any integer solution with the value $(k_2,\tau_2)$ has the objective value at most as good as the corresponding symmetry strategy with $(k_2,\tau_2)$. This significantly helps the algorithm use a \textit{branch and bound} technique inside itself.

The pseudocode for solving Model 2 is similar to Model 1. In this model, any pair of $(k,\tau)$ that is able to satisfy the $k \times P_{time} + m \times T_{time} \leq n \times \tau$ (i.e. for at most $p=1$), coupled with its corresponding symmetry solution of $G,D$ satisfying $Pr(infection|(k,\tau, G,D)) \leq \alpha \times background~risk$, is a potential feasible solution and can be explored by the heuristic search approach to find the minimum $p$ value.
For a given pair $(k,\tau)$ and its corresponding symmetry solution, the heuristic search first finds the $O(2^k)$ possible integer neighbors of $D$. That means, for any $d_i \in D$, there are two possibilities to rounding, $floor (d_i)=\lfloor d_i \rfloor$ and $ceil(d_i)=\lceil d_i \rceil$, except $d_k$, which always to be assumed $d_k=\tau$.  In terms of time complexity, this is possible in a reasonable time for the real-world sizes of the problem.  For example, if the residents are tested once a week in the worst case ($Max_{\tau}=7$), the maximum possible solutions for rounding the testing days never exceeds 64. Otherwise, in theory, and for large-scale instances (i.e., for $\tau \geq k>15$), we may choose just a subset of possible neighbors.

Note that, since we assume one round of tests per day, the number of groups is less or equal to $\tau$. Indeed, $floor (d_i)$ and $ceil(d_i)$ are the integers with at most one unit distance from $d_i$ in the integer space. 
However, we can extend the idea for more than one unit distance, but based on our brute force searches for practical cases, doing at least one test every two weeks per resident, the optimal solution is never more than one distance unit away. In the next step, the algorithm explores possible integer combinations of $G$. Similar to the rounding set $D$, the algorithm searches over the size of the groups in the integer space but, in this case, for more than one unit. 

For example, for input parameters $m=30, Max_g=20, n=5, P_{time}=180~minutes, T_{time}= 15~minutes, p=10\%, \kappa=0.5 m$, and $\beta=0.1$ with the 7-day incidence per 1000 individuals, the optimal solutions of Model 1 is $\tau=5, k=4$ and $D=\{d_1=1,d_2=2,d_3=4,d_4=5\}$ and $G=\{|g_1|=6,|g_2|=10,|g_3|=8,|g_3|=6\}$ with the optimal objective value, \textit{expected detection time} of 1.2133 days. In this optimal solution, the size of groups is $10-\frac{30}{4}=2.5$ units far away from the symmetry solution's group size. The optimal solution of Model 2 for exactly the same input parameters with the infection level $\alpha=0.5$, is the symmetry strategy $\tau=3, k=3$ and $D=\{d_1=1,d_2=2,d_3=3\}$ and $G=\{|g_1|=10,|g_2|=10,|g_3|=10\}$ with the optimal value of \textit{expected detection time} 0.8658 days, \textit{risk of infection} 0.0007, and the objective value $p=13.75\%$. In the next section, we will provide more results of the algorithm for different settings of the input parameters.

\section{Simulation Results and Discussion}
\label{Simulation}

\begin{figure}[H]
\begin{subfigure}{.5\textwidth}
  \centering
  \includegraphics[width=0.9 \linewidth]{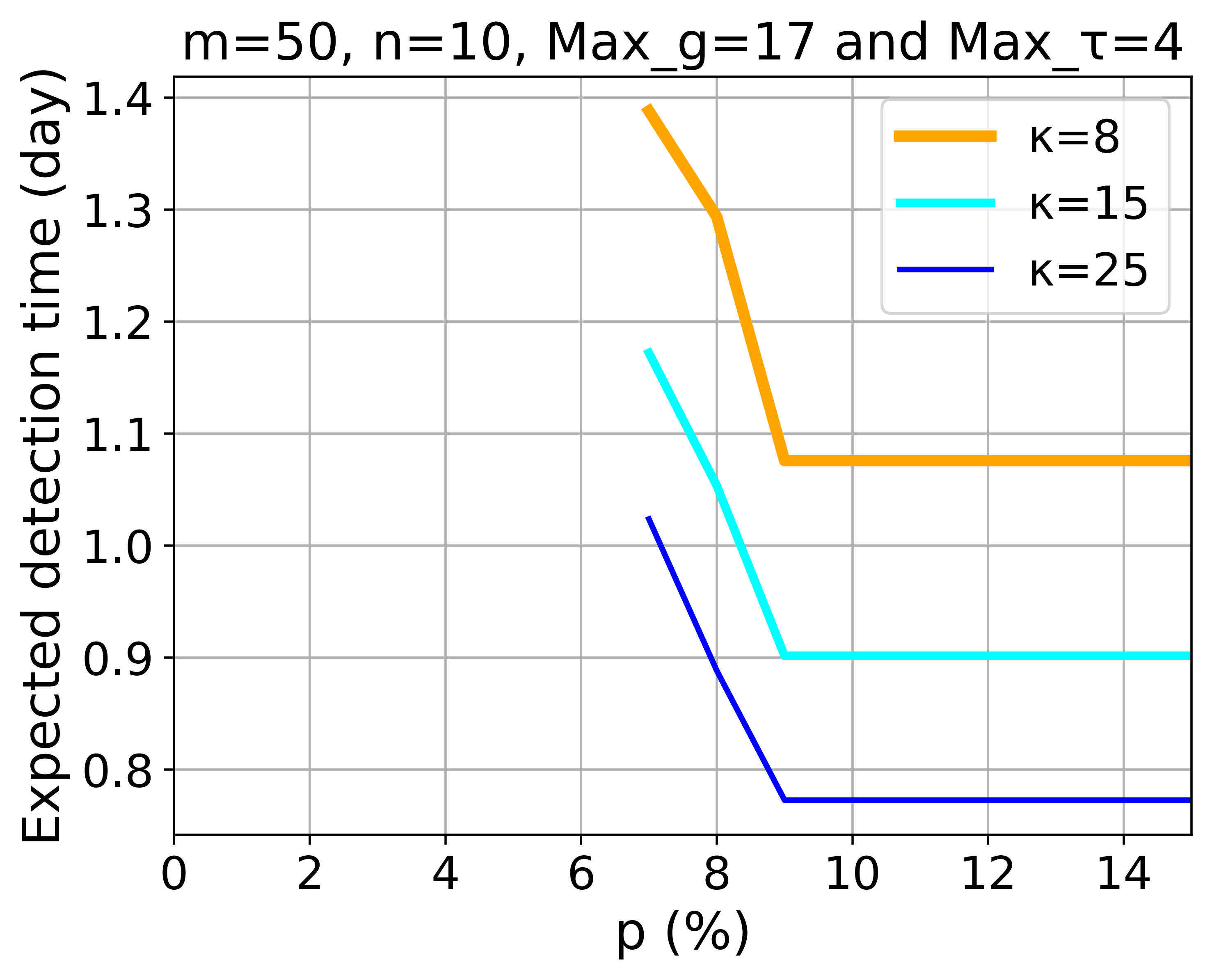}  
\end{subfigure}
\begin{subfigure}{.5\textwidth}
  \centering
  \includegraphics[width=0.9 \linewidth]{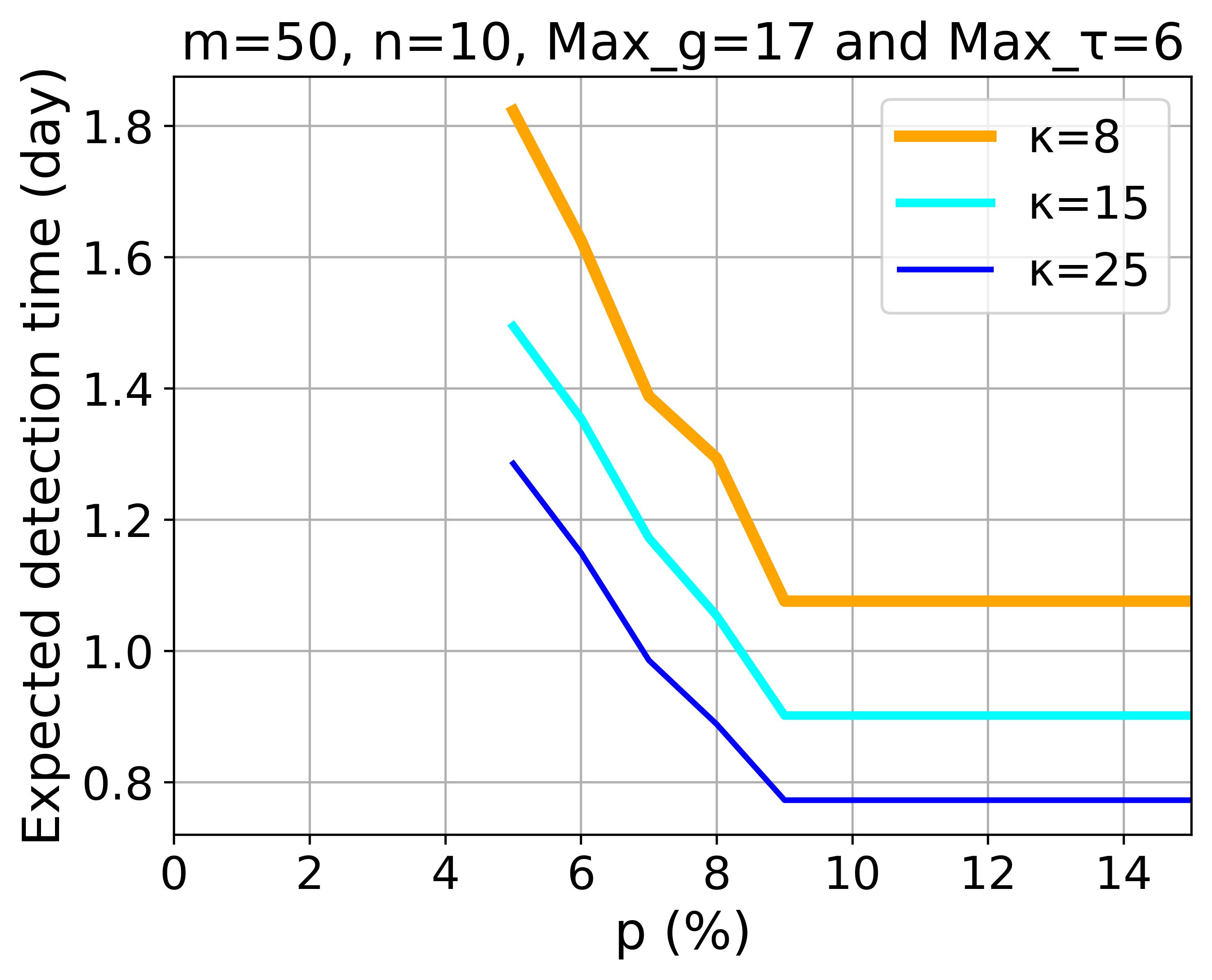}  
\end{subfigure}
\begin{subfigure}{.5\textwidth}
  \centering
   \includegraphics[width=0.9 \linewidth]{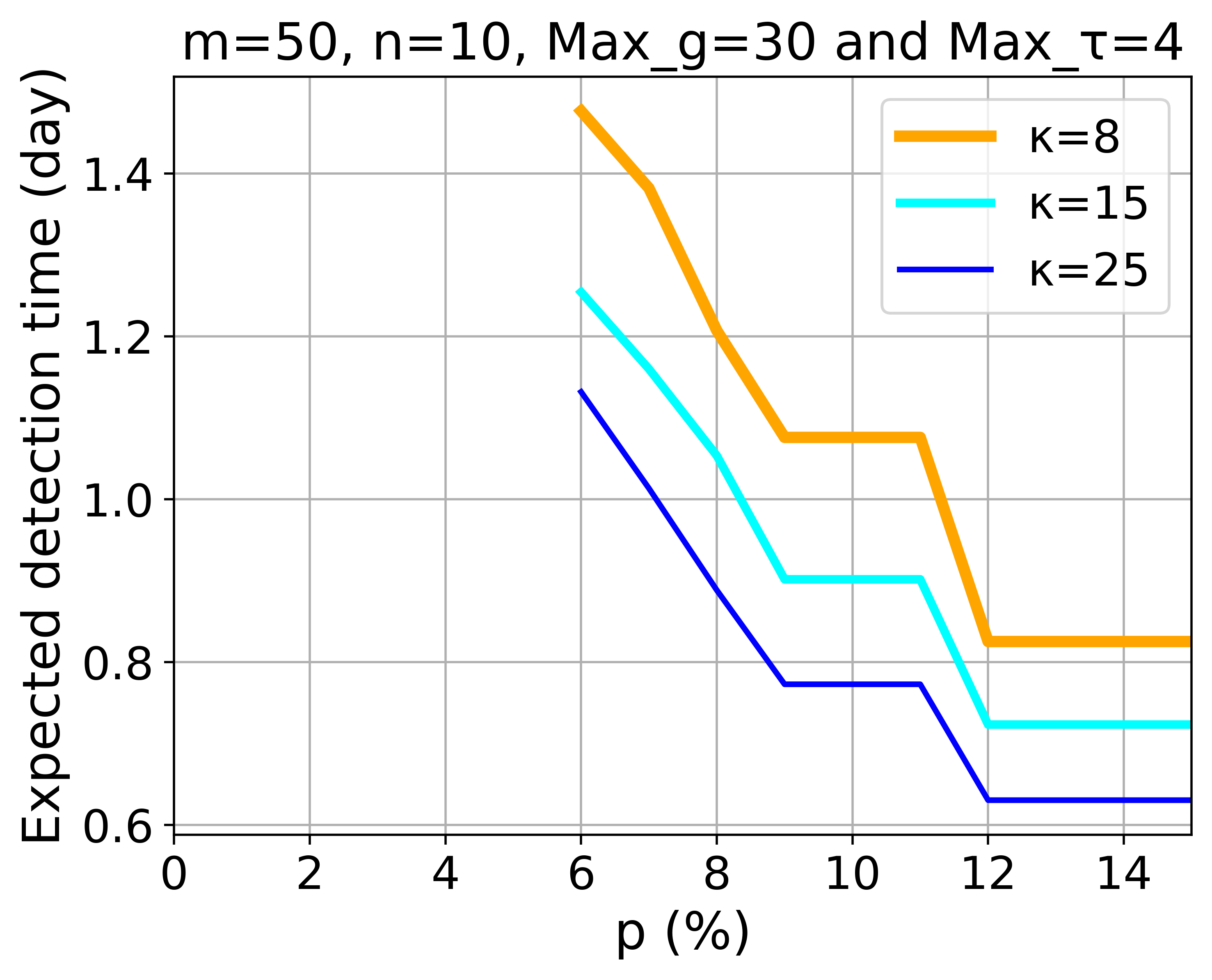}  
\end{subfigure}
\begin{subfigure}{.5\textwidth}
  \centering
   \includegraphics[width=0.9 \linewidth]{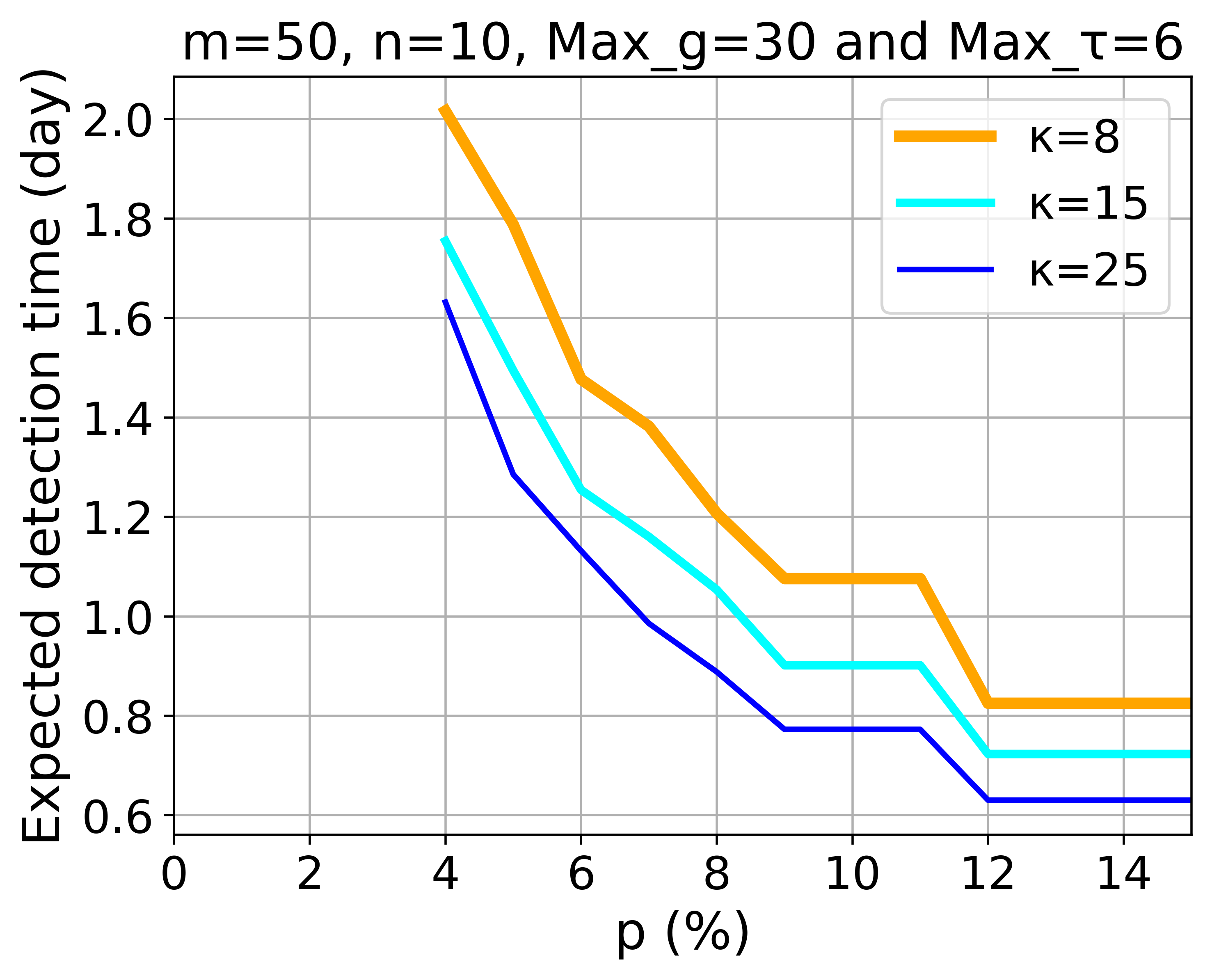}  
\end{subfigure}
\begin{subfigure}{.5\textwidth}
  \centering
  \includegraphics[width=0.9 \linewidth]{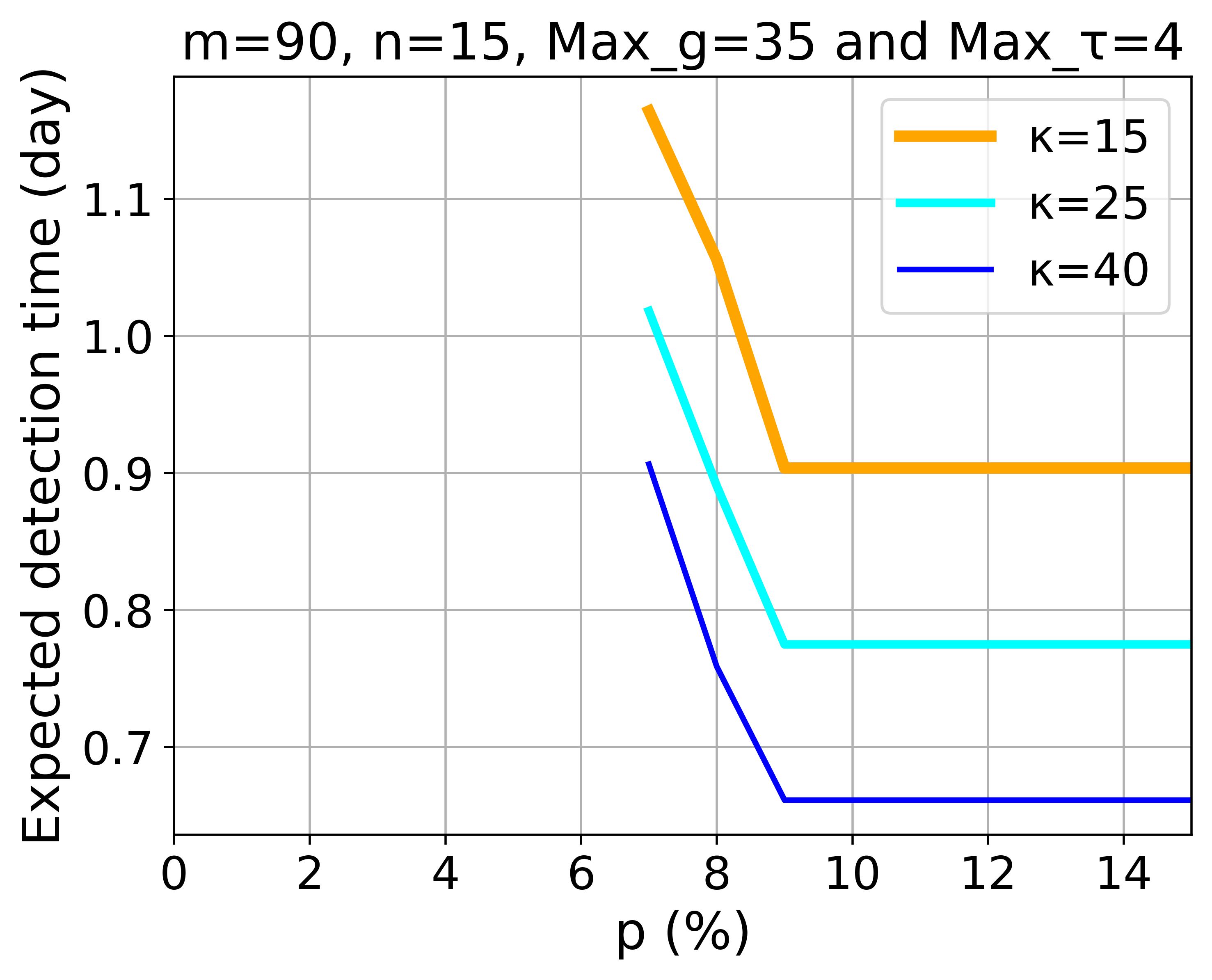}  
\end{subfigure}
\begin{subfigure}{.5\textwidth}
  \centering
  \includegraphics[width=0.9 \linewidth]{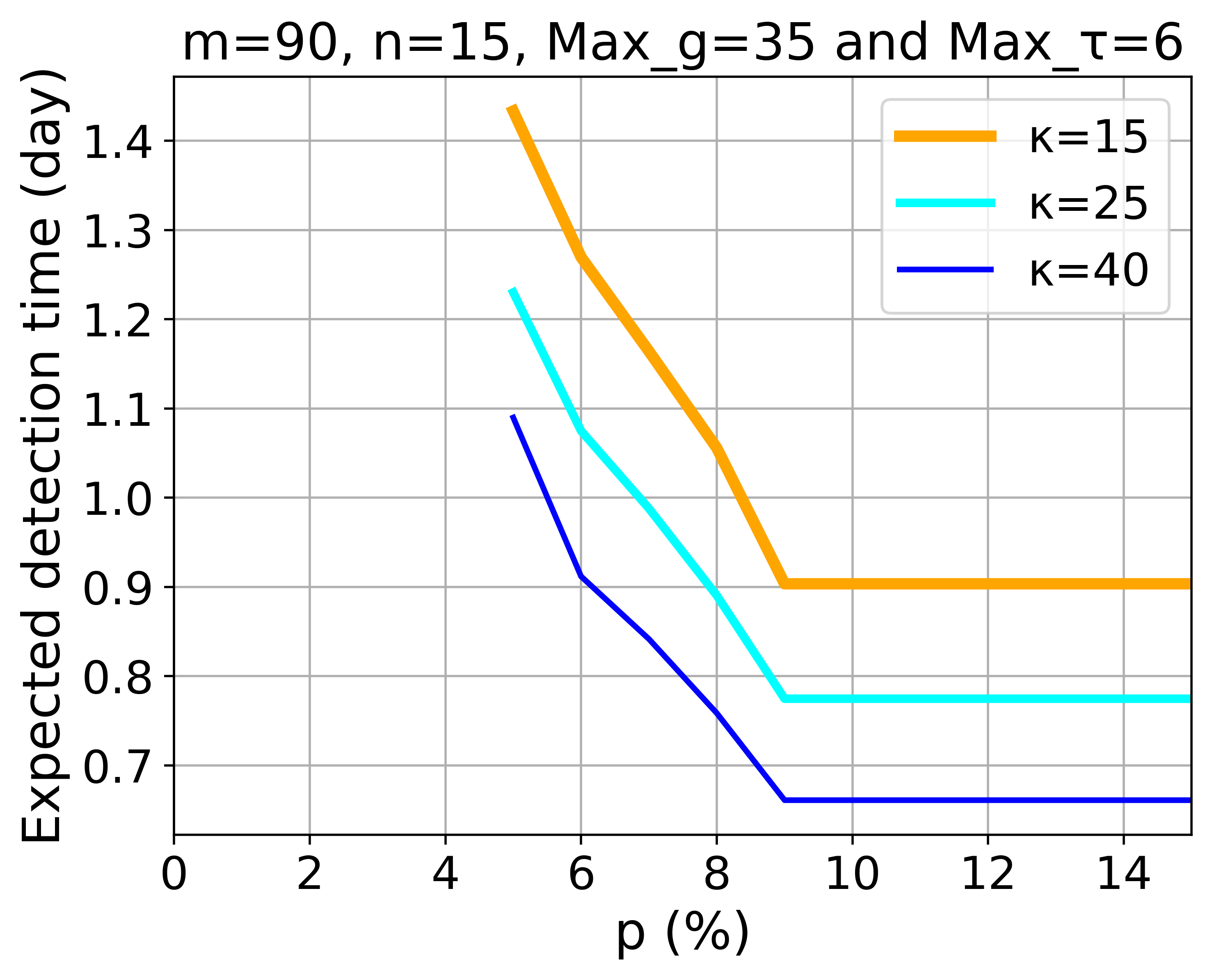}  
\end{subfigure}
\begin{subfigure}{.5\textwidth}
  \centering
   \includegraphics[width=0.9 \linewidth]{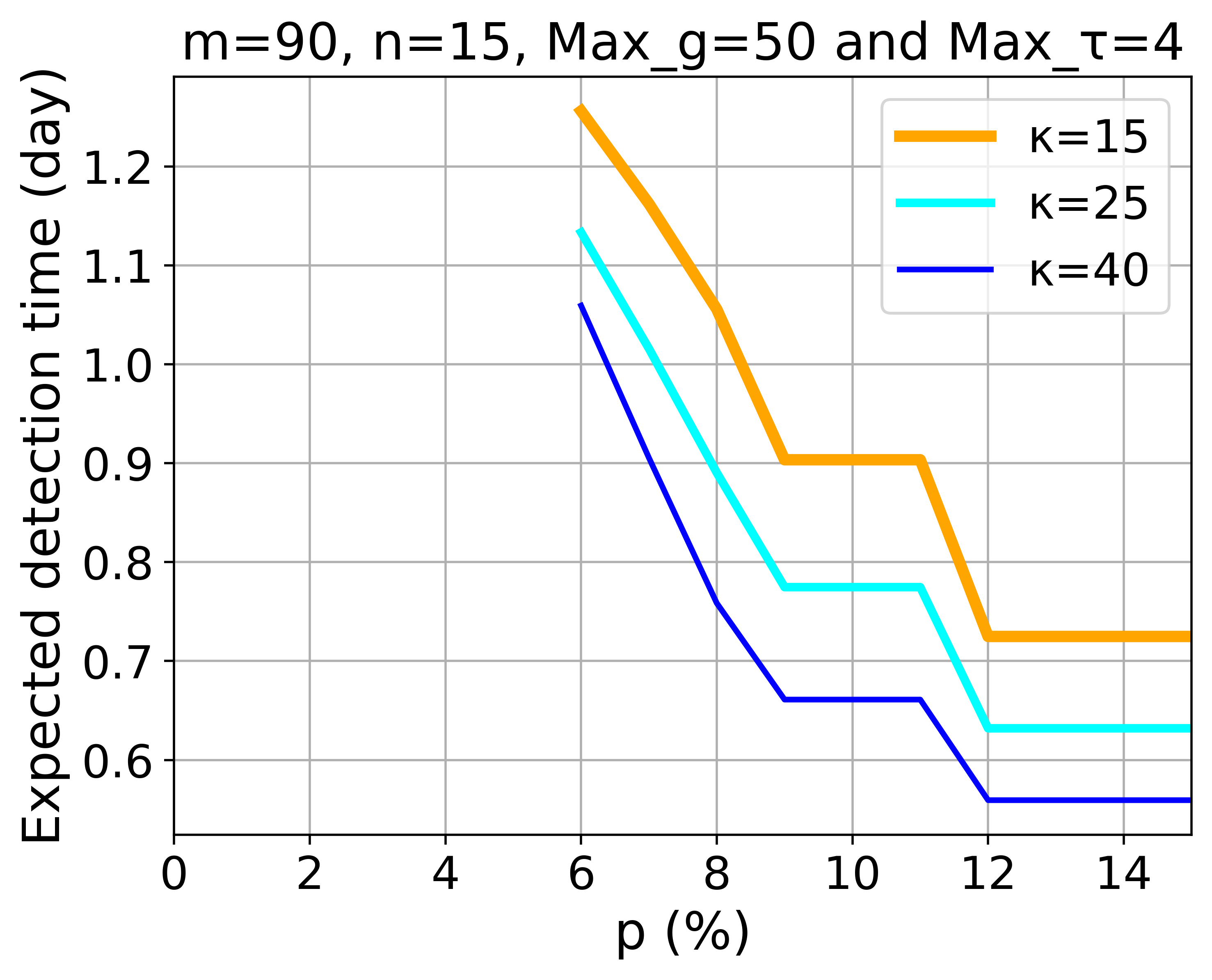}  
\end{subfigure}
\begin{subfigure}{.5\textwidth}
  \centering
   \includegraphics[width=0.9 \linewidth]{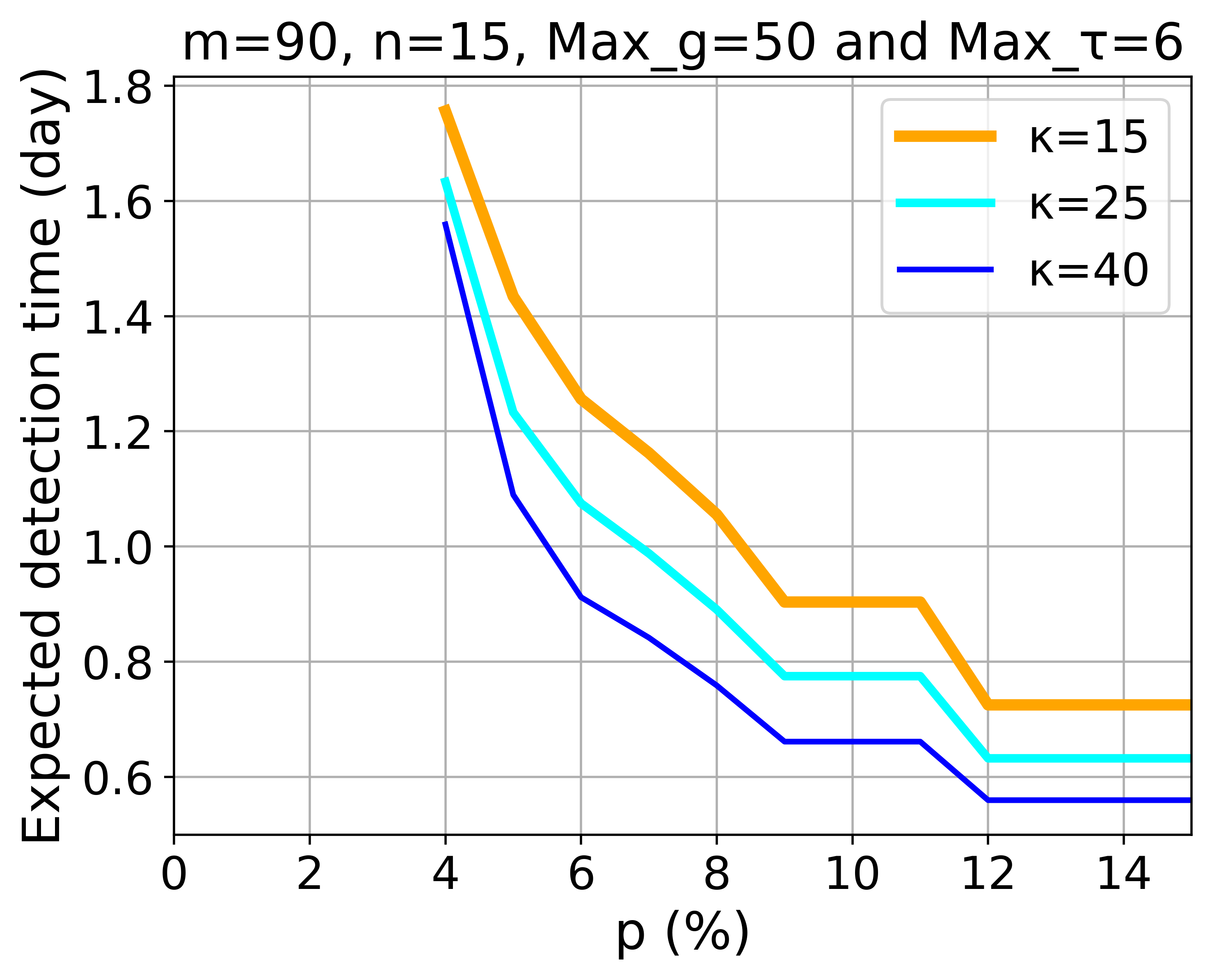}  
\end{subfigure}
\caption{Results of Model 1 for different settings of input parameters. Each subfigure shows the obtained expected detection time for different values of $p \in [0,15]$ for three different number of contacts per resident per day.}
\label{Pareto_Model1}
\end{figure}

\begin{figure}[H]
\begin{subfigure}{.5\textwidth}
  \centering
  \includegraphics[width=0.9 \linewidth]{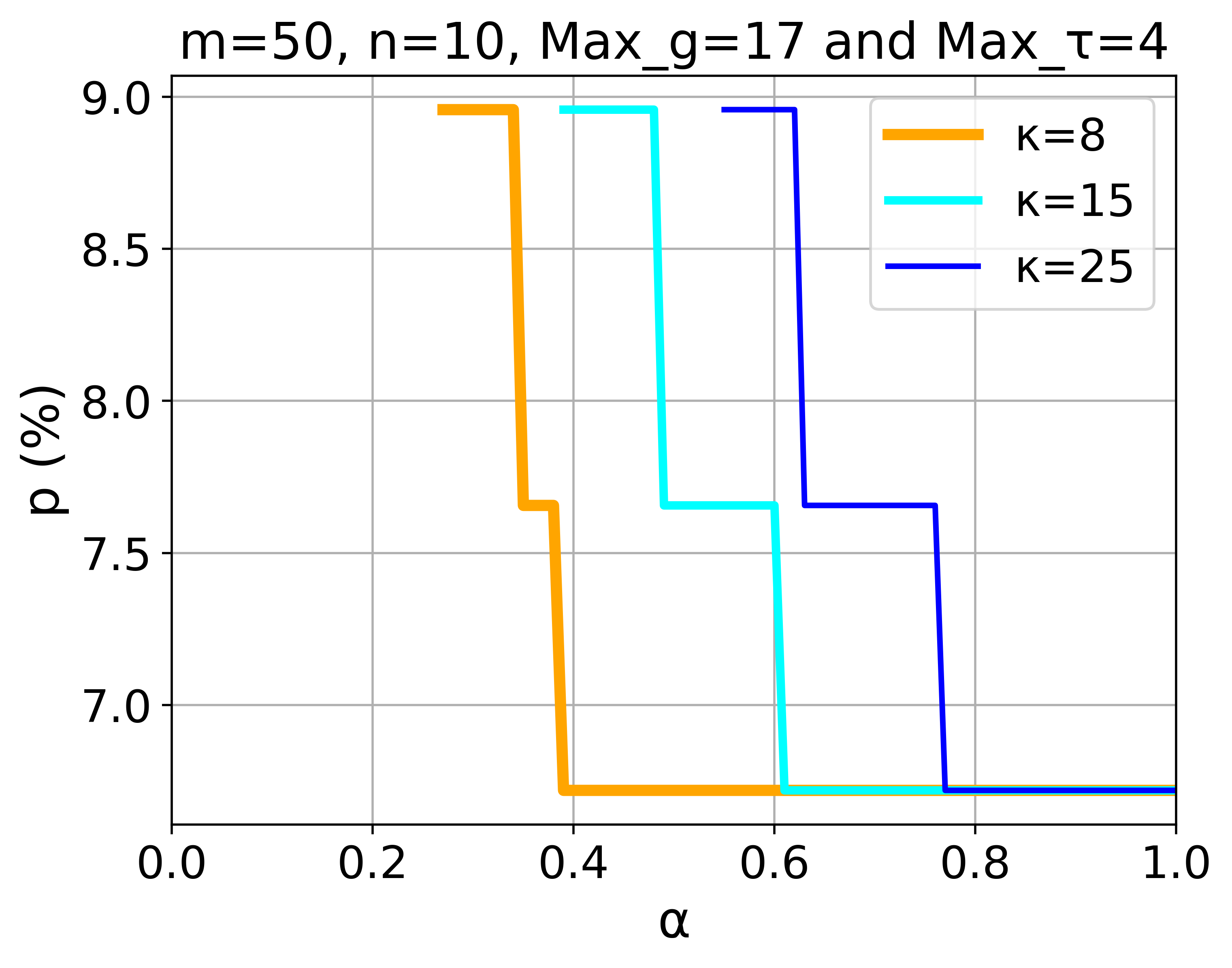}  
\end{subfigure}
\begin{subfigure}{.5\textwidth}
  \centering
  \includegraphics[width=0.9 \linewidth]{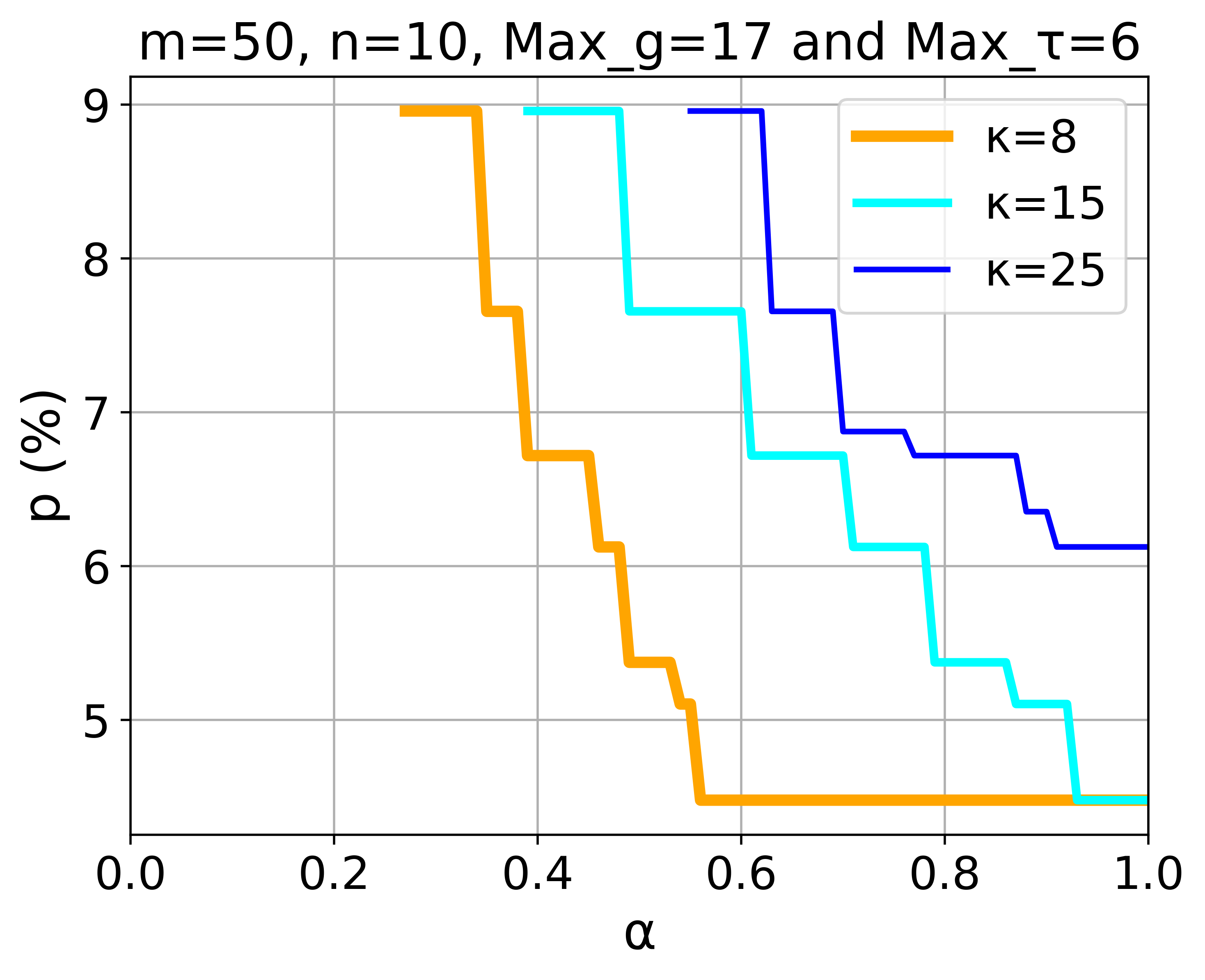}  
\end{subfigure}
\begin{subfigure}{.5\textwidth}
  \centering
   \includegraphics[width=0.9 \linewidth]{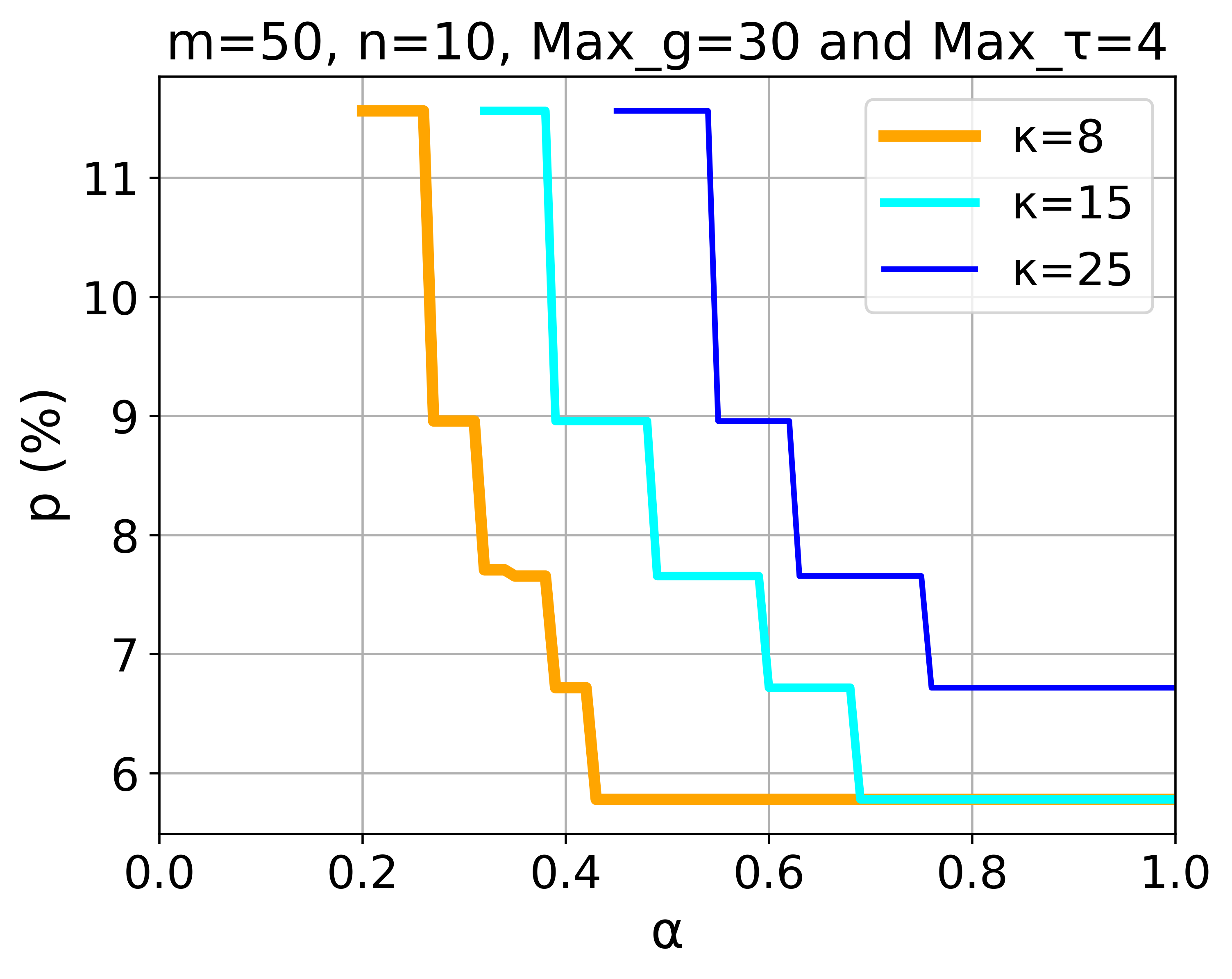}  
\end{subfigure}
\begin{subfigure}{.5\textwidth}
  \centering
   \includegraphics[width=0.9 \linewidth]{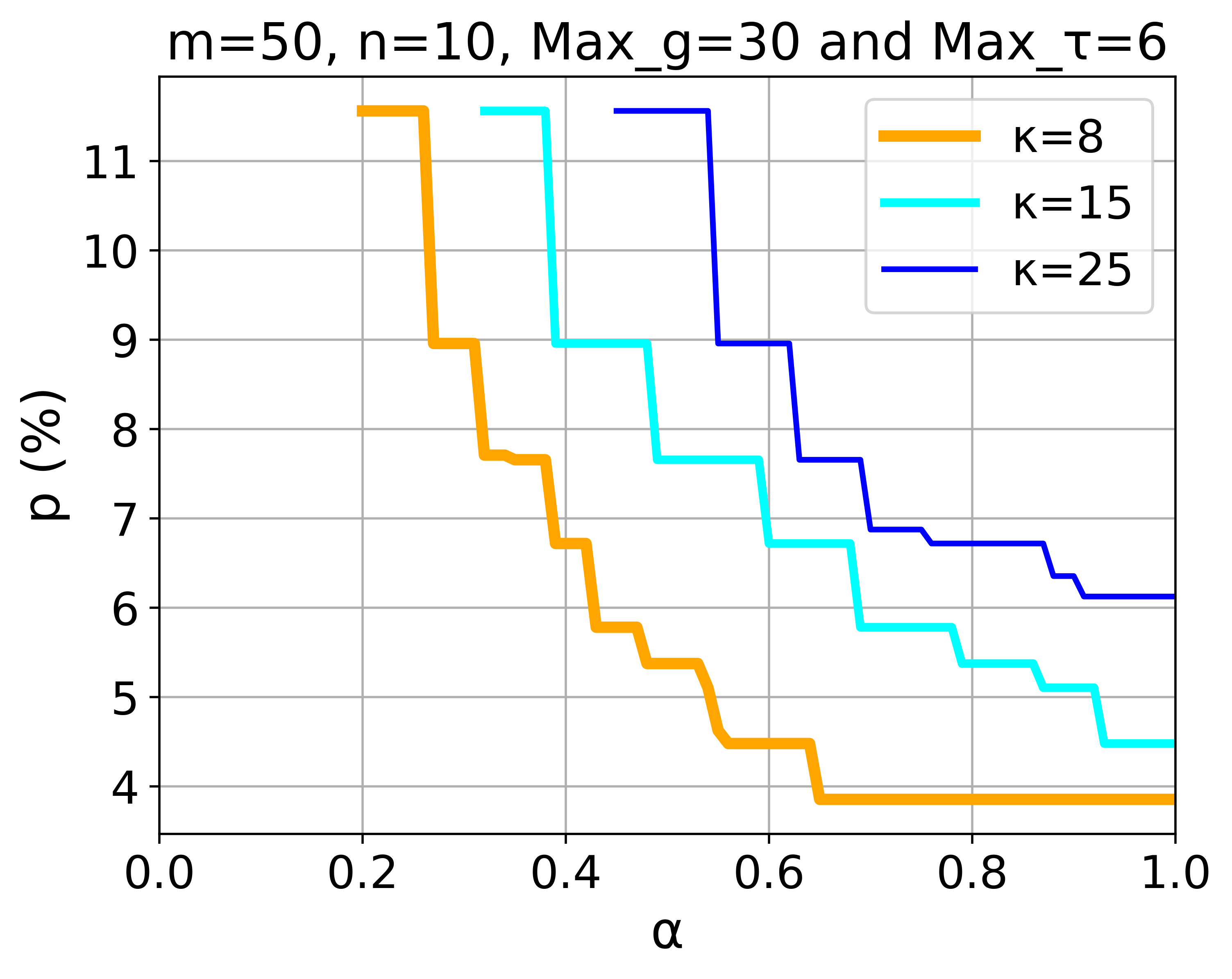}  
\end{subfigure}

\begin{subfigure}{.5\textwidth}
  \centering
  \includegraphics[width=0.9 \linewidth]{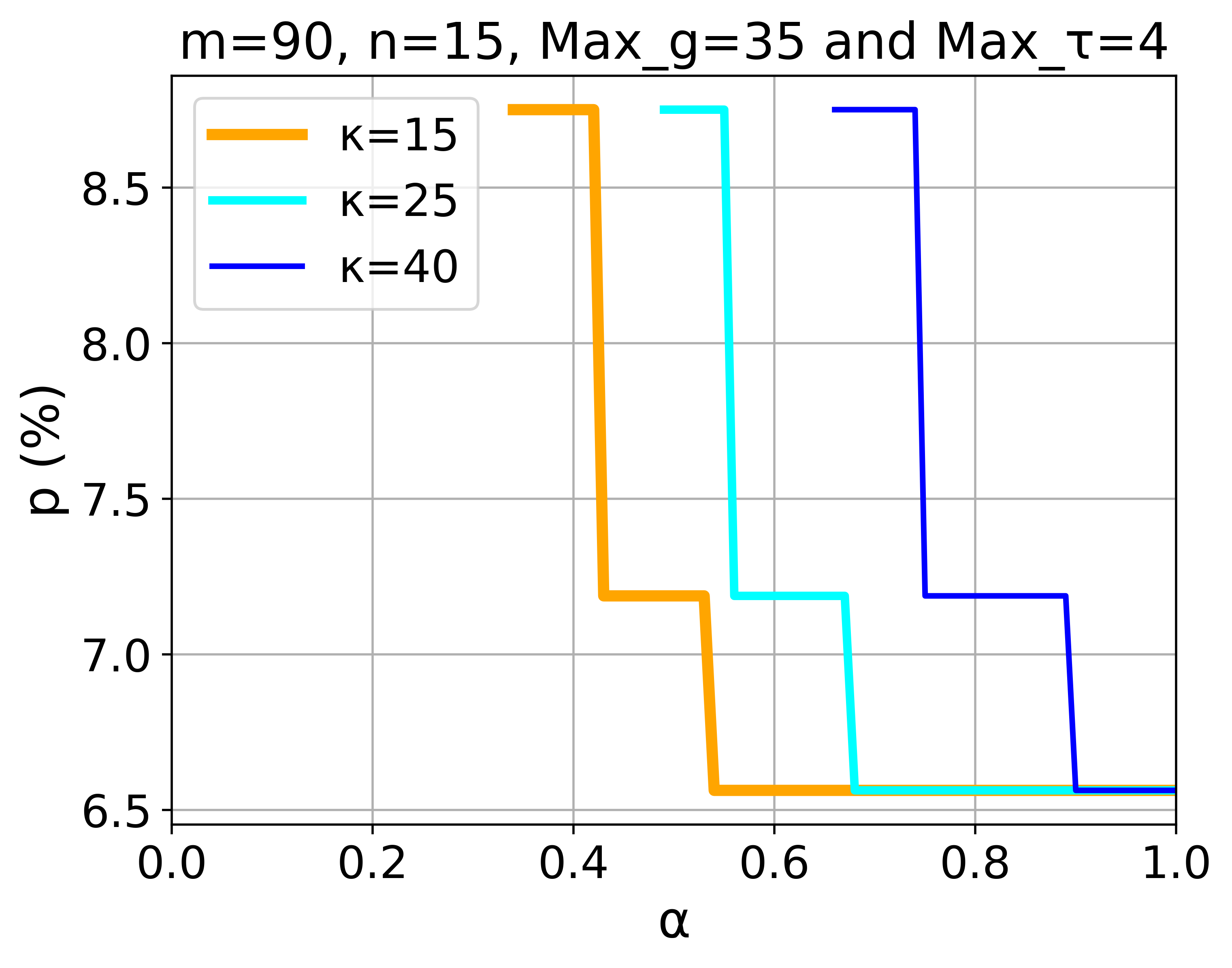}  
\end{subfigure}
\begin{subfigure}{.5\textwidth}
  \centering
  \includegraphics[width=0.9 \linewidth]{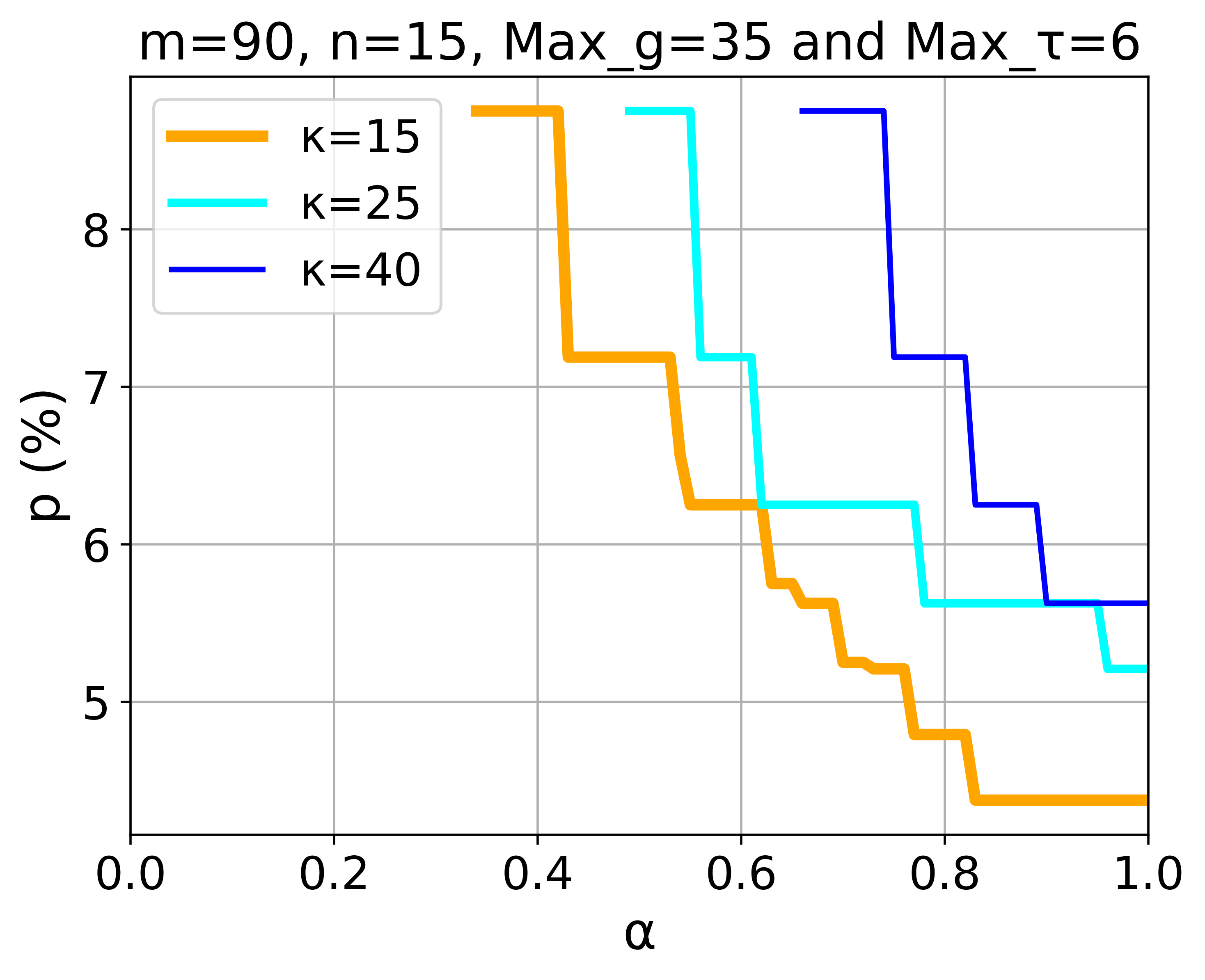}  
\end{subfigure}
\begin{subfigure}{.5\textwidth}
  \centering
   \includegraphics[width=0.9 \linewidth]{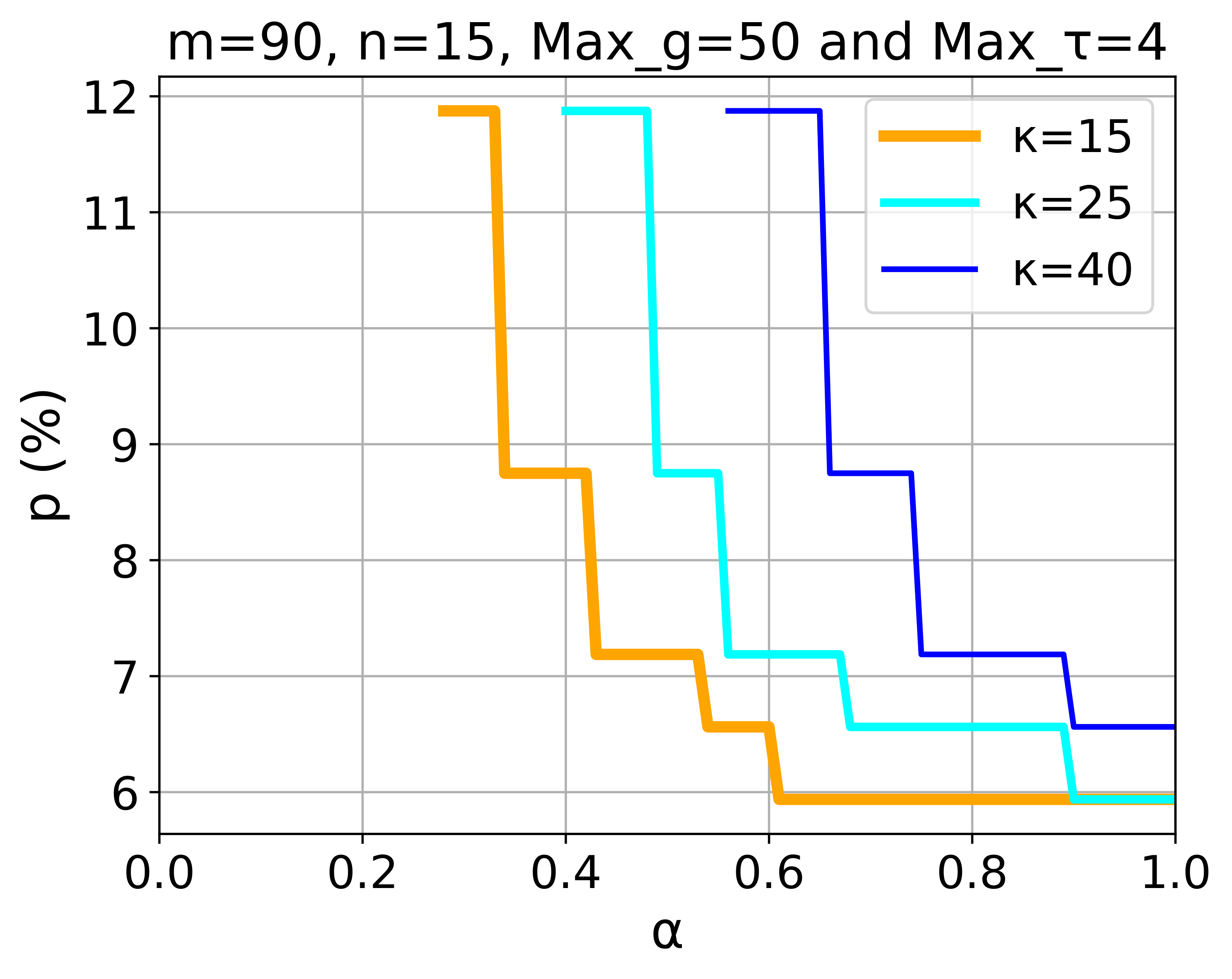}  
\end{subfigure}
\begin{subfigure}{.5\textwidth}
  \centering
   \includegraphics[width=0.9 \linewidth]{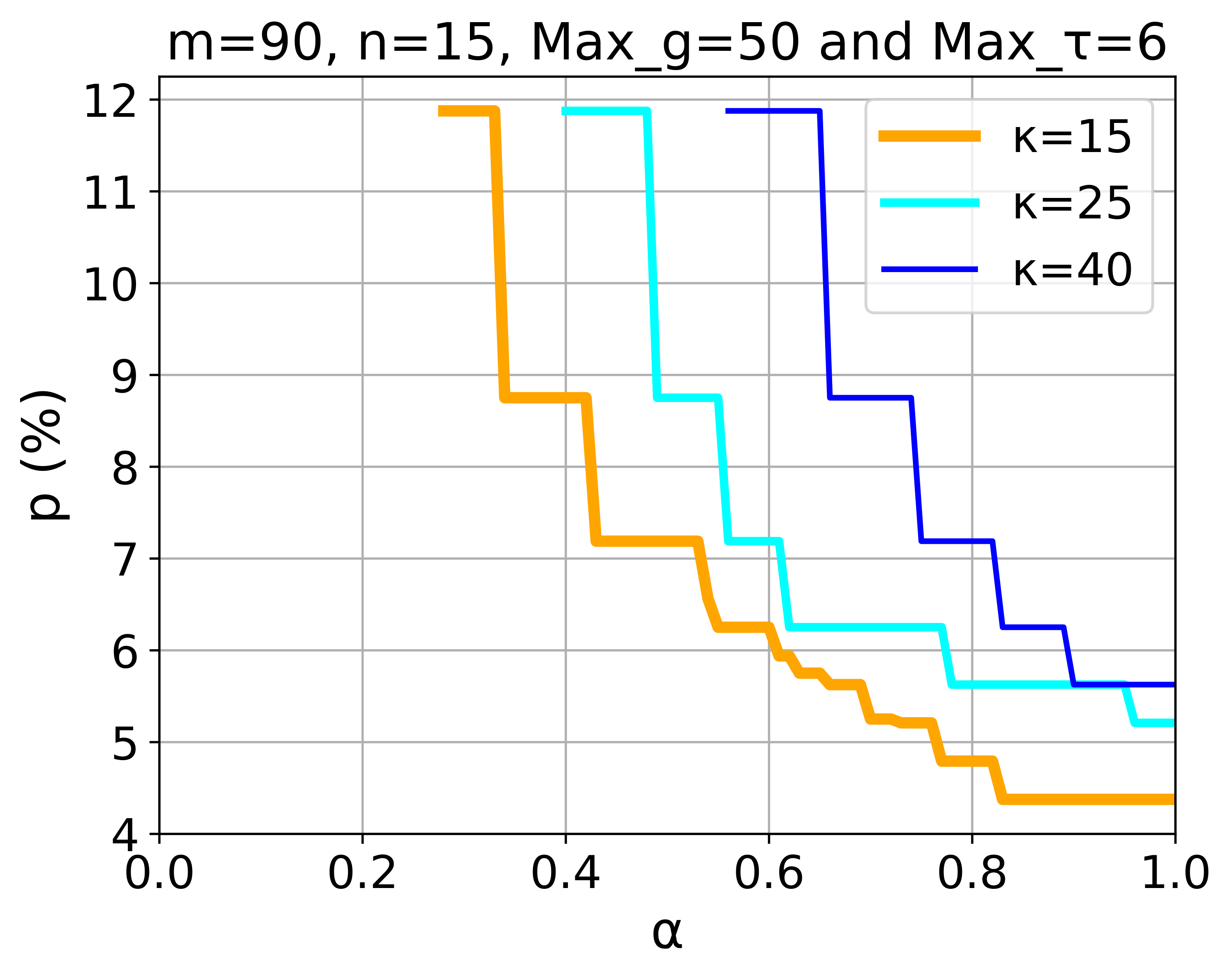}  
\end{subfigure}

\caption{Results of Model 2 for different settings of input parameters. Each subfigure shows the obtained portion of staff's time for different values of $\alpha \in [0,1]$ for three different number of contacts per resident per day.}
\label{Pareto_Model2}
\end{figure}

This section is organized in three parts and shows the results of the proposed models and algorithm to find the optimal testing strategy in an RH. The first and second parts shows the optimal testing strategy $(\kappa,\tau,G,D)$ for 48 different settings of input parameters. The results are shown in Table \ref{Table1} and Table \ref{Table2} in the Appendix, and a subset of them is discussed in Table \ref{Table1_subset} and Table \ref{Table2_subset}. The third part (Figure \ref{Pareto_Model1} and Figure \ref{Pareto_Model2}) presents a sensitivity analysis between the main decision parameter and the objective values, That is, the tradeoff between the expected detection time and the staff workload, $p$, in Model 1, and the tradeoff between the expected detection time and the staff workload, $p$, in Model 1, and the tradeoff between $p$ and the risk of infection in RH, $\alpha$, in Model 2.
The code of the algorithm is implemented in the programming language Python 3.7 and runs on a standard PC (\textit{Intel}(\textit{R}) \textit{Core}(\textit{TM}) \textit{i}7 and 32G \textit{RAM}). We remark that the proposed solution approach can solve real-world scale problems in a short time. For a RH with $m=90$ residents, $n=15$ staff, and a testing frequency of at least once a week, the run time is about 7 seconds. 

The experiments are divided into two parts. The first shows the resulting testing strategies, including $k,\tau, G, D$, and the objective values for realistically sized combinations of the inputs. 
In the second part, we perform a sensitivity analysis over the parameter $p$ in Model 1 and parameter $\alpha$ in Model 2. So, we illustrate Pareto optimal solutions of the proposed models based on these values.

For obtaining the optimal testing strategy, we set the input parameters to realistic values as follows:

\begin{itemize}
    \item $(m,n) \in \{(50,10),(90,15)\}$
    \item $Max_\tau \in \{4,7\}$
    \item $Max_g \in \{\lceil \frac{m}{2} \rceil,\lceil \frac{m}{3} \rceil \}$
    \item For Model 1, $p \in \{5\%, 10\%,20\%\}$
    \item For Model 2, $\alpha \in \{30\%, 50\%,75\%\}$
    \item The number of contacts per resident per day, (\textit{i}) $\kappa=0.15 \times m +5$, and (\textit{ii}) $\kappa=0.3 \times m +5$
    \item The preparation time, $P_{time}=180$ minutes and the testing time, $T_{time}=15$ minutes
    \item The probability of disease transmission per contact, $\beta=0.1$. We choose this value as a probable pessimistic case from a possible range of values reported in previous studies~\cite{lelieveld2020model,yang2021sars, world2022enhancing, lyngse2021sars} regarding the first variants of COVID-19 (\textit{Alpha}, \textit{Delta}, and \textit{Omicron}, which is more transmissible than the previous ones).

    \item \textit{The background risk} is calculated considering 7-day incidences of COVID-19 infections as 600 individuals per 100,000 population (Reference data of weekly incidences in Saxony, Germany in the period of October 15 to December 15, 2021.)
\end{itemize}

%

We ran Model 1 and Model 2 for 48 different combinations of the input parameters. Furthermore, we assumed that a probable infection could arrive at an RH by the $n$ staff members or by visitors. We considered one visitor on average per resident every two weeks. Given a 7-day incidence value of 600, the probability of infection per individual per week is simply computed as $\frac{600}{100,000}$. Consequently, the probability of an infection arriving at the RH is $1-(1-\frac{600}{100000})^{n+\frac{m}{14}}$. 

The full results tables for Model 1 and Model 2 are contained in the Appendix. For all the reported results, we ran the proposed algorithm 5 times per input and returned the best obtained solution. Out of $5 \times 48 = 240$ independent runs, the algorithm succeeded in reaching the global optimal solution of the problem in 236 runs, and in the remaining 4 runs, it obtained local optimal solutions very close to the global ones. The optimal solutions are computed by a brute force algorithm. 

\paragraph{\textbf{Results Model 1}}
Table \ref{Table1} in the Appendix summarize the full results of Model 1. Here, we focus on a subset (see Table \ref{Table1_subset}) of runs that highlight important properties of the model. The \textit{Run 1} - \textit{Run 8} represent the results for the case in which the decision-maker allocates at most 5\% of staff working load to testing. We observe that for the \textit{Run 1} - \textit{Run 4}, there is no feasible solution when it is preferred to test the residents once every 4 days. 
However, when it is extended to test once a week, there is always a feasible testing strategy. If the staff can test 30 residents per day and each of the residents has 9 contacts per day (see \textit{Run 5}), the optimal solution is the non-symmetry solution $\tau=5$, $k=2$, $G=\{28,22\}$ and $G=\{2,5\}$ with minimum expected detection time equal to 1.7365 days. In contrast, if the staff can test 22 residents per day (see \textit{Run 7}), the objective value is increased to 1.7587 days, and the optimal solution is $\tau=6$, $k=3$, $G=\{16,17,17\}$ and $G=\{2,4,6\}$.

Another comparison can be made between \textit{Run 5} and \textit{Run 6}, where the only difference between them is the number of contacts, $\kappa$. When $\kappa=9$, the optimal expected detection time is 1.7365 days, while for the contact number $\kappa=17$ it is 1.4355 days. As described, in the first case, the optimal solution is choosing a test interval $\tau=5$, partitioning the residents into two groups with sizes 28 and 22, and testing them on days 2 and 5. While for the second one, $\tau=6$ is the optimal test interval, and it is better to partition into three groups with sizes 16, 17 and 17, and test them on days 2, 4, and 6. Note that the second case is the closest solution to the corresponding symmetry solution, while the first one is relatively far from its corresponding symmetry solution. In fact, the symmetry solution is $\{2.5,5\}$ with the same group size equal to 25. So, since $\{2.5,5\}$ is infeasible, the algorithm rounded it to $D=\{2,5\}$, and in proportion to such testing days, it changed the group size to 28 and 22 residents to achieve the minimum possible expected detection time. To draw an analogy between this obtained optimal solution by the algorithm and the solution $D=\{2,5\}$ and $G=\{25,25\}$, the objective value for the optimal one is 1.7365 days, while for the latter one is 1.74198 days.

Another result is also found in \textit{Run 29} which is for $m=90$ residents with a test interval of at most 7 days. The proposed algorithm chooses 6 as the test interval, partitions the residents into 4 groups with sizes $G=\{25,20,25,20\}$, and tests them on days $D=\{1,3,4,6\}$. For this setting, the resulting expected detection time is 1.4332 days. However, partitioning like $\{22,23,22,23\}$ and testing days $\{2,3,5,6\}$, which is a more even and uniform distribution, results in a higher value of the expected detection time (1.4343 days).

\begin{table}[]
\caption{Summary of discussed results of 48 runs for Model 1 (see the extended Table \ref{Table1} in the appendix). Columns $C_5$, $C_6$ and $C_{12}$ correspond to the parameters $Max_{\tau}$, $Max_g$ and the expected detection time, respectively. The runs with no feasible solution are shown by blank cells.}
\label{Table1_subset}
\tiny
\begin{tabular}{|l|l|l|l|l|l|l|l|l|l|l|l|}
\hline
       & \textbf{m}  & \textbf{n}  & \textbf{p}    & \textbf{$C_5$} & \textbf{$C_6$} & \textbf{$\kappa$} & \textbf{$k$} & \textbf{$\tau$} & \textbf{$G$}                       & \textbf{$D$}                 & \textbf{$C_{12}$} \\ \hline

Run 1  & 50 & 10 & 0.05 & 4        & 30     & 9        &   &     &                       &                 &           \\ \hline
Run 2  & 50 & 10 & 0.05 & 4        & 30     & 17       &   &     &                       &                 &           \\ \hline
Run 3  & 50 & 10 & 0.05 & 4        & 22     & 9        &   &     &                       &                 &           \\ \hline
Run 4  & 50 & 10 & 0.05 & 4        & 22     & 17       &   &     &                       &                 &           \\ \hline
Run 5  & 50 & 10 & 0.05 & 7        & 30     & 9        & 2 & 5   & \{28,22\}             & \{2,5\}         & 1.7365    \\ \hline
Run 6  & 50 & 10 & 0.05 & 7        & 30     & 17       & 3 & 6   & \{16,17,17\}          & \{2,4,6\}       & 1.4355    \\ \hline
Run 7  & 50 & 10 & 0.05 & 7        & 22     & 9        & 3 & 6   & \{16,17,17\}          & \{2,4,6\}       & 1.7587    \\ \hline
Run 29 & 90 & 15 & 0.05 & 7        & 30     & 15       & 4 & 6   & \{25,20,25,20\}       & \{1,3,4,6\}     & 1.4332    \\ \hline
\end{tabular}
\end{table}

\paragraph{\textbf{Results Model 2}}
The results of Model 2 are presented in Table \ref{Table2} in the Appendix and a subset of them are presented in Table \ref{Table2_subset}. In this model, the parameter $\alpha$ plays an important role in defining the feasible and infeasible space. So, as can be seen in the table when a small value is selected for $\alpha$ (i.e., $\alpha =0.3$), in most cases there is no feasible solution for the given input settings. However, for bigger values such as $\alpha =0.75$, there is always at least one feasible solution. The portion of the staff's workload, which can be allocated to the testing process of the resident, changes from 3.84 to 8.96. The maximum value is related to minimum $\alpha=0.3$. That means that if the decision-maker wants to achieve a low infection risk level, such as 0.3 of background risk, it is necessary to allocate at least $8.96\%$ of the staff's time to testing (e.g., \textit{Run 3}). In contrast, the value $\alpha=0.75$ can result in a portion like $3.84\%$ for some settings of the inputs (See \textit{ Run 23}).

This is because this model is not directly focused on the objective of the expected detection time and as soon as a feasible solution (i.e., a testing strategy that satisfies Constraint \eqref{M2.10}, $Pr(infection|(k,\tau, G, D)) \leq \alpha$), is obtained, it tries to minimize the portion of staff time allocated to the testing process. So, some particular values of the risk level $\alpha$ may change the boundary of the feasible and infeasible solution space. For example, for a setting of input parameters such as $m=50, n=10, \kappa=9, Max_g=30$ and $Max_{\tau}=5$, the optimal solution for $\alpha=0.60$ is $k=4$, $\tau=4$, $G=\{12,13,13,13\}$ and $D=\{1,2,3,5\}$, with the objective value $p=5.104\%$, which is a symmetry solution. While, if we set the risk level to $\alpha=0.61$, the optimal solution will be the non-symmetry strategy $k=2$, $\tau=5$, $G=\{26,24\}$ and $D=\{2,5\}$, with the objective value $p=4.625\%$. A similar case happens for Model 1 for particular values of the parameter $p$. So, in the next round of simulations, we illustrate the trade-off between parameter $p$ as the main constraint of Model 1 and the expected detection time, as well as the trade-off between parameter $\alpha$ as the main constraint of Model 2 and the portion of staff time which must be allocated to the testing process.

Overall, the set of solutions shown in Table \ref{Table2} can help decision-makers to choose the desired optimal strategy, considering the available resources in the RH and the local incidence situation. For example, the results of the last four runs of Table \ref{Table1} are the same optimal strategy either $G=\{30,30,30\}$ and the testing days $D=\{1,2,3\}$, or $G=\{18,18,18,18,18\}$ and the testing days $D=\{1,2,3,4,5\}$, with objective values between 0.7381 and 1.1904 days. These results can be achieved by at most 20\% of the staff's workload. The corresponding results in Table \ref{Table2}, show solutions for $\alpha=0.75$, while all the solutions follow the pattern $G=\{18,18,18,18,18\}$ and the testing days $D=\{1,2,3,4,5\}$. As it is clear, among these solutions, the ones in which the resident has 15 contacts per day result in a better objective value of $5.21\%$.
These findings can be significantly valuable for managers of RHs to establish new regulations for the contacts among residents, including restrictions for visitors and isolation measures, updating them efficiently over time.


\begin{table}[]
\caption{Summary of discussed results of 48 runs for Model 2 (see the extended Table \ref{Table2} in the appendix). Columns $C_5$, $C_6$ and $C_{12}$ corresponds to the parameters $Max_{\tau}$, $Max_g$, and the obtained minimum portion of the time for the testing process of the residents which should be allocated by the staffs, respectively. 
}
\label{Table2_subset}
\tiny
\begin{tabular}{|l|l|l|l|l|l|l|l|l|l|l|l|}
\hline
       & \textbf{$m$}  & \textbf{$n$}  & \textbf{$\alpha$}    & \textbf{$C_5$} & \textbf{$C_6$} & \textbf{$\kappa$} & \textbf{$k$} & \textbf{$\tau$} & \textbf{$G$} & \textbf{$D$} & \textbf{$C_{12}$} \\ \hline
Run 3  & 50 & 10 & 0.3   & 4        & 22     & 9        & 3 & 3   & \{16,17,17\}       & \{1,2,3\}     & 8.96    \\ \hline
Run 21 & 50 & 10 & 0.75  & 7        & 30     & 9        & 3 & 7   & \{16,17,17\}       & \{2,4,7\}     & 3.84    \\ \hline
Run 23 & 50 & 10 & 0.75  & 7        & 22     & 9        & 3 & 7   & \{16,17,17\}       & \{2,4,7\}     & 3.84    \\ \hline
Run 41 & 90 & 15 & 0.75  & 4        & 30     & 15       & 3 & 4   & \{30,30,30\}       & \{1,2,4\}     & 6.56    \\ \hline
Run 45 & 90 & 15 & 0.75  & 7        & 30     & 15       & 5 & 6   & \{18,18,18,18,18\} & \{1,2,3,4,6\} & 5.21    \\ \hline
Run 46 & 90 & 15 & 0.75  & 7        & 30     & 29       & 5 & 5   & \{18,18,18,18,18\} & \{1,2,3,4,5\} & 6.25    \\ \hline
Run 47 & 90 & 15 & 0.75  & 7        & 22     & 15       & 5 & 6   & \{18,18,18,18,18\} & \{1,2,3,4,6\} & 5.21    \\ \hline
Run 48 & 90 & 15 & 0.75  & 7        & 22     & 29       & 5 & 5   & \{18,18,18,18,18\} & \{1,2,3,4,5\} & 6.25    \\ \hline
\end{tabular}
\end{table}

\paragraph{\textbf{Sensitivity analysis}}
In the second part of the experiments, we performed a sensitivity analysis to evaluate the trade-off between the expected detection time and the staff workload, $p$. For Model 2 we analyze the trade-off between the level of acceptable risk in the facility, $\alpha$, and the staff workload. 

Figure \ref{Pareto_Model1} and Figure \ref{Pareto_Model2} show the sensitivity analysis for Model 1 and Model 2, respectively. We evaluated the case of $m=50$ residents with $\kappa= 8, 15$ and $250$, and the case $m=90$ residents with $\kappa= 15, 24$ and $40$. The value of parameters $m, n, Max_g$ and $Max_{\tau}$ for each run are reported on the top of each subfigure.

Figure \ref{Pareto_Model1}, compares different values of $p$, the portion of staff's time which can be allocated to the testing process (horizontal axis), and the expected detection time (vertical axis). The diagrams are shown for $p \leq 15$ (i.e., the results for $p > 15$ are the same as for $p = 15$). 
We observe that there is no feasible solution for small values of $p$. We also note that there is a maximum value of $p$ to minimize the expected detection time of a probable infection for a given set of input parameters, and it is independent of the contact number. So, a decision-maker may choose such a value of $p$ as the best solution to minimize the detection time without considering the number of contacts among the residents. 
For example, in the first subfigure (the case $m=50, n=10, Max_g=17$ and $Max_{\tau}=4$, $p=9$)
, there is no gain for the RH to increase $p$ more than 9\% so that 91\% of staff's time can be allocated to their caring tasks. Also, if the retirement home can test more residents in one day (or say for $Max_g$), then it is possible to minimize the expected detection time to reach the ideal expected value of 0.5, say in 12 hours. Note that we assumed one round of testing per day is possible, so the ideal detection time can never be less than 0.5 days (i.e., 12 hours). Having a large group size, i.e., $Max_g \geq \frac{m}{2}$, and a large value of $p$ will result in achieving such an ideal detection time. 

Overall, the reported results in Figure \ref{Pareto_Model1} display the trade-off between the risk of infection for residents, $\alpha$, and the portion of time which the staff can allocate for the testing process, $p$. Indeed, since $1-p$ is the portion of staff's time allocated to caring tasks, we may interpret both of these results from the residents' point of view: the risk of infection and comfort level of the residents. 
So, for a given configuration of the input parameters, by decreasing the risk of infection of the residents, their comfort level decreases. In this result, the number of contacts plays an important role, and decreasing it helps decrease the risk of infection and, consequently, decreases $p$. For example, in the last subfigure, $p$ varies from 12 to almost 4; for a risk value of $\alpha=0.4$, there is no solution for high contact numbers $\kappa=40$, $p=12$ for the medium contact numbers $\kappa=25$, and $p \approx 9$ for the low contact numbers $\kappa=15$. Similarly to the case shown in Figure \ref{Pareto_Model1}, among all 8 subfigures in Figure \ref{Pareto_Model2}, the minimum objective value $p \approx 4$ is reached for the cases that the retirement home has the capability of testing at least half of the resident in one day ($Max_g \geq \frac{m}{2}$).
As expected, by increasing the rate of contact between the residents in a RH, a probable infection will also transmit quickly, and consequently, the time to detect it by the testing process will decrease. However, this does not mean the number of infected individuals will decrease. For example, in the first subfigure in Figure \ref{Pareto_Model1}, the expected detection time for $\kappa=25$ is less than for $\kappa=8$ for any value $p$. 

Figure \ref{Pareto_Model2} illustrates the trade-off between $\alpha$ and $p$; the level of risk of infection and the staff workload and for different contact numbers. The results are shown for $\alpha \in [0,1]$. As we mentioned before, note that for small values of $\alpha$, there is no feasible solution. For example, in the first subfigure ($m=50$ residents, $n=10$ staff, $Max_g=17$ and $Max_{\tau}=4$), there is no feasible solution for $\alpha < 0.27$ when the residents have $\kappa=8$ contact per day on average. Moreover, if they have $\kappa=15$ contacts, there is no feasible solution for $\alpha < 0.38$. Finally, for $\kappa=25$, the minimum $\alpha$ which results in a feasible solution is $0.54$. Note that, for the corresponding setting of the input in Figure \ref{Pareto_Model2}, for the case $\kappa=25$, the rate of infection for any particular value $\alpha$ is higher than for the case $\kappa=8$.

Moreover, in the first subfigure, it is clear that when $\alpha$ changes, the portion of staff time only changes 2.5 units (from 9\% to 6.5\%). That means allocating 9\% of staff workload to the testing process is more reasonable because it significantly reduces the risk of infection. Finally, it is worth mentioning that the diagrams are decreasing but not monotone. Indeed, some critical values (the bending points of the diagrams) of $\alpha$ resulted in a better portion of the testing workload. That means, it is not to be expected that by increasing any small value in the portion of the testing workload of the staff, a better level of risk value is obtained. Thus, the decision-maker can consider the critical values of $\alpha$ to allocate an efficient portion of staff workload to the testing process.

\section{Conclusion and future work}
\label{Conclusions}

In this paper, we developed two novel MINLP models to compute the optimal testing strategy for residents in retirement homes during the COVID-19 pandemic. The models aimed to minimize the risk of infection and therefore prevent the spread of it, considering the trade-off between the portion of the staff's workload allocated to the testing process and the frequency of tests for the residents. Because of the residents' high risk of mortality given infection, any step to shorten detection time is vital and could result in saving lives.

Model 1 derived an optimal testing strategy by minimizing the expected detection time of the virus in the facility while considering the maximum threshold on the staff workload. In other words, the manager obtains the testing strategy which minimizes the infection risk for any given threshold on the workload. 
Model 2 determined an optimal testing strategy when the manager specifies a threshold as the maximum tolerable risk of infection inside the RH. As a reference point, we considered the risk of infection in the neighborhood of the retirement home.

We proposed a practical and efficient approach to solve the models. We discussed the properties of the optimal testing strategy and observed that the symmetry solution is optimal for the models in the continuous search space. Based on this observation, we presented a global local search heuristic algorithm to find an optimal testing strategy in the integer search space.
We verified the models and the algorithm by testing several realistically sized instances of the problem. In most instances, the optimal solutions are symmetry strategies in the integer search space, but the algorithm found both solutions with and without the symmetry property.

The results of the experiments revealed that, in a retirement home, an optimal testing strategy depends not only on the local incidence level, but also on the rate of contact among the residents and the maximum size of a group of residents that the staff can test in one batch, especially when the optimal solution is asymmetric. Thus, decision-makers can use the proposed models to build a set of optimal solutions that better suit their expectations and available resources. These findings are novel because, to date, there is no other alternative to assess the effectiveness of testing schedule strategies in retirement homes to minimize the resident's risk of infection during a pandemic. 

We identified several future research directions. Firstly, in this study, we assumed homogeneity of the residents in terms of the probability of infection. However, this is a general case that is not always realistic in practice as the residents differ by age, co-morbidity, and vaccination status, which may alter their probability of infection. So, considering heterogeneity in the probability of infection is an interesting aspect to study. Secondly, we considered two possibilities for an infection arriving at the retirement home--via the staff members and visitors. However, it is possible to consider additional sources of infection, such as physiotherapists or doctors that regularly visit the residents. 
Thirdly, in the modeling phase, we considered the \textit{average} value of the input parameters (i.e., number of contacts and the probability of transmitting the virus), so the proposed models and the algorithms work well for the instances of the problem with low variability. However, for instances with high variability, the computed expected detection time might not work properly. So, the proposed model can be extended to cover such cases.
Finally, future work can consider more realistic scenarios, incorporating different shifts and testing strategies for the staff and considering a contact network that reflects the social interactions between the staff and residents in a retirement home.


\section*{Acknowledgements}
This work was partially funded by the Where2Test project, which is financed by SMWK with tax funds on the basis of the budget approved by the Saxon State Parliament. This work was also partially funded by the Center of Advanced Systems Understanding (CASUS) which is financed by Germany’s Federal Ministry of Education and Research (BMBF) and by the Saxon Ministry for Science, Culture and Tourism (SMWK) with tax funds on the basis of the budget approved by the Saxon State Parliament.
Further, the authors wish to thank the managers of the retirement homes Diakonisches Werk im Kirchenbezirk Löbau-Zittau GmbH, Saxony, Germany, for their collaborating in providing data and insights used to solve the problem presented in this paper.

\bibliography{bibfile}

\section*{Appendix}
\label{Appendix}

\textbf{Theorem} \ref{thm_eq}.
In the problem of finding an optimal testing strategy for a retirement home assuming homogeneity of residents, the objective of minimizing the detection time is equivalent to minimizing the number of infections after an infection arrives at the retirement home.

\begin{proof}
\label{thm_eq_proof}
Suppose at day $d_0$ an infection arrives at the facility, and after $d$ days, it is detected (no matter the source of infection or any other resident is detected). Assuming a complete contact network and homogeneous rate of interactions among the residents. The probability of infection for any resident is obtained by Equation \ref{Prob_Infection} (See Subsection \ref{Computingprobability}), the function $P_I(m,\kappa,\beta,d)$ which is a strictly increasing function to $d$. On the other hand, the expected number of infected residents can be computed by $1 + (m-1) \times P_I(m,\kappa,\beta,d)$. Consequently, minimizing the detection time $d$ results in minimizing the expected number of infections.~~~~~~$\blacksquare$
\end{proof}

\textbf{Theorem} \ref{thm_limit}
The optimal testing strategy of the real space is a symmetry strategy for the cases $\beta \longrightarrow 0$, $\beta \longrightarrow 1$, $\kappa \longrightarrow 0$, or $\kappa \longrightarrow +\infty$.

\begin{proof}
\label{thm_limit_proof}
For a testing strategy $(k,\tau,G,D)$, where $G=\{g_1,g_2,\dots,g_k\}$, and $D=\{d_1,d_2,\dots,d_k\}$, let define $l_1=d_1$ and $l_i=d_i-d_{i-1}$, for $i=2,3,\dots,k$ as the distance (time) between the testing day of each group and its preceding group.

In the cases of $\beta \longrightarrow 1$ or $\kappa \longrightarrow +\infty$, the probability of infection transmits from the source of infection to the other residents will go to zero. Precisely, we can write the probability of infection, Eq. \ref{Prob_Infection},  as follows

\begin{dmath}
\lim_{\beta \to 1} P_I(m,\kappa,\beta,d)=\lim_{\kappa \to +\infty} P_I(m,\kappa,\beta,d)=
1-(1-P_I(m,\kappa,\beta,d-1)) \times \left[(1-\beta)^{\kappa \frac{1}{m-1}} \times (1-P_I(m,\kappa,\beta,d-1)\beta)^{\kappa (1-\frac{1}{m-1})} \right]= 1-(1-P_I(m,\kappa,\beta,d-1)) \times (0)=1.
\label{lim1_Prob_Infection}
\end{dmath}

Thus, the expected detection time, Eq. \ref{Exp_time}, can be computed as

\begin{equation}
Expected~Detection~ Time~ of~(k,\tau, G,D)= \sum_{i=1}^{k} \left( \frac{l_i}{\tau} \sum_{t_0=1}^{l_i} t_0 \right)= \frac{1}{2\tau} \sum_{i=1}^{k} l_i(l_i+1).
\end{equation}

Regarding the fact $\sum_{i=1}^{k} l_i= \tau$, the above equation can be written as
\begin{equation}
Expected~Detection~ Time~ of~(k,\tau, G,D)= \frac{1}{2\tau} \sum_{i=1}^{k} l_i(l_i+1) = \frac{1}{2\tau} \left(  \tau+ \sum_{i=1}^{k} l_i^2 \right).
\end{equation}

Therefore, to minimize the detection time, we need to minimize $\sum_{i=1}^{k} l_i^2$ under the condition $\sum_{i=1}^{k} l_i= \tau$. Now, for a symmetry testing strategy, we have $l_1=l_2=\dots=l_k=\frac{\tau}{k}$. By contradiction, suppose this property does not hold in some optimal testing strategy $S$. Thus, there exists some pair of intervals $l_i$ and $l_j$ in $S$ such that $l_i \neq l_j$. Without changing the value of the other $\tau-2$ intervals, we set $l_i=l_j=\frac{l_i+l_j}{2}$ and show that a better testing strategy can be obtained, which means, $S$ is not optimal. Contradiction. Also, it is clear that $l_i^2+l_j^2>2(\frac{l_i+l_j}{2})^2$, because

\begin{equation}
l_i^2+l_j^2>2(\frac{l_i+l_j}{2})^2 \Longleftrightarrow 2l_i^2+2l_j^2>l_i^2+l_j^2+2l_il_j\Longleftrightarrow (l_i-l_j)^2>0.
\end{equation}

Thus, the proof is complete for the cases $\beta \longrightarrow 1$ and $\kappa \longrightarrow +\infty$. Now, let consider the cases of $\beta \longrightarrow 0$ or $\kappa \longrightarrow 0$. In these cases, the probability of infection transmits from the source of infection to the other residents will approach zero. Precisely, we can write the probability of infection, Eq. \ref{Prob_Infection},  as follows: 

\begin{dmath}
\lim_{\beta \to 0} P_I(m,\kappa,\beta,d)=\lim_{\kappa \to 0} P_I(m,\kappa,\beta,d)=
1-(1-P_I(m,\kappa,\beta,d-1)) \times \left[(1-\beta)^{\kappa \frac{1}{m-1}} \times (1-P_I(m,\kappa,\beta,d-1)\beta)^{\kappa (1-\frac{1}{m-1})} \right]=1-(1-P_I(m,\kappa,\beta,d-1)) \times (1)=P_I(m,\kappa,\beta,d-1).
\label{lim0_Prob_Infection}
\end{dmath}

and since $P_I(m,\kappa,\beta,0)=0$, we can conclude
\ref{Prob_Infection},  as follows
\begin{equation}
\lim_{\beta \to 0} P_I(m,\kappa,\beta,d)=\lim_{\kappa \to 0} P_I(m,\kappa,\beta,d)=0.
\end{equation}

This indeed means that for detecting the infection, we need to detect the source of the infection. Therefore, the expected detection time, Eq. \ref{Exp_time}, can be computed as

\begin{dmath}
Expected~Detection~ Time~ of~(k,\tau, G,D)= \frac{1}{\tau} \sum_{t_0=0}^{\tau-1} \mathbb{E} (k,\tau, G,D,t_0)=\frac{1}{\tau} \sum_{t_0=0}^{\tau-1}\frac{\tau}{2}=  \frac{\tau}{2}.
\end{dmath}

That means no matter which $k$, $G$, or $D$ are chosen, the only important parameter to minimize the expected detection time is minimizing test interval $\tau$. So, again, a Symmetry strategy is an optimal solution for the problem.~~~~$\blacksquare$
\end{proof}

\begin{algorithm} 
\caption{Proposed algorithm for solving the problem in Model 1}
\textbf{Input:} Retirement home's parameters ($m,n,P_{time}, T_{time}, Max_{\tau}, Max_g$, $p$)\\
\textbf{Output:} A testing strategy to minimize the expected detection time\\ 
\begin{algorithmic}[1]
\STATE $opt\_exp\_time \leftarrow +\infty$
\STATE $opt\_strategy \leftarrow \varnothing$
\FOR{$\tau=1$ \textbf{to} $Max_{\tau}$}
\FOR{$k=\lceil \frac{m}{Max_g} \rceil$ \textbf{to} $\tau$}
\IF{$k \times P_{time} + m \times T_{time} \leq p \times n \times \tau$}
\STATE $G=\{|g_1|=\frac{\tau}{k},|g_2|=\frac{\tau}{k},\dots,|g_k|=\frac{\tau}{k}\}$
\STATE $D=\{d_1=1\frac{\tau}{k},d_2=2\frac{\tau}{k},\dots,d_k=k\frac{\tau}{k}\}$
\STATE $S \leftarrow (k,\tau,G,D)$
\STATE $exp\_time \leftarrow$ expected detection time of $S$ (Eq. \ref{Exp_time})
\IF{$exp\_time < opt\_exp\_time$}
\IF{$S$ is not an integer solution}
\STATE $(S,exp\_time) \leftarrow \textbf{Heuristic\_Search}(S,m,n,P_{time}, T_{time}, Max_{\tau}, Max_g, p)$
\ENDIF
\IF{$exp\_time < opt\_exp\_time$}
\STATE $opt\_exp\_time \leftarrow exp\_time$
\STATE $opt\_strategy \leftarrow S$
\ENDIF
\ENDIF
\ENDIF
\ENDFOR
\ENDFOR
\IF{$opt\_strategy == \varnothing$}
\STATE \textbf{print} (\textit{There is no feasible strategy})
\ELSE
\STATE \textbf{return} $S$ and $opt\_exp\_time$
\ENDIF 
\end{algorithmic}
\label{AlgorithmFirstModel}
\end{algorithm}

\begin{algorithm}
\caption{Heuristic\_Search}
\textbf{Input:} A non-integer Strategy $S=(k,\tau,G,D)$ and the retirement home's parameters ($m,n,P_{time}, T_{time}, Max_{\tau}, Max_g$, $p$)\\
\textbf{Output:} An integer testing strategy in the neighbor set of $S$ with the minimum expected detection time
\begin{algorithmic}[1]
\STATE $best\_exp\_time \leftarrow +\infty$
\STATE $best\_stratey \leftarrow \varnothing$
\STATE $Max\_interation \leftarrow k \times \tau \times m$
\STATE $All\_Days \leftarrow$ All integer neighbors of $D$ with less than one unit distance
\FOR {$ any D' \in All\_Days$}
\STATE $G' \leftarrow$ The nearest integer neighbors of $G$ such that $\sum_{i=1}^{k} |g_i| =m$
\STATE $exp\_time_{G'} \leftarrow$ The expected detection time of $(k,\tau,G',D')$ (Eq. \ref{Exp_time})
\STATE $T_0 \leftarrow exp\_time_{G'}$
\FOR{$i=1$ \textbf{to} $Max\_interation$}
\IF{$exp\_time_{G'} < best\_exp\_time$}
\STATE $best\_exp\_time \leftarrow exp\_time_{G'}$
\STATE $best\_strategy \leftarrow (k,\tau,G',D')$
\ENDIF
\STATE $G'' \leftarrow$ A random integer neighbor of $G'$ with less than one unit distance from $G'$ such that $\sum_{i=1}^{k} |g_i| =m$
\STATE $exp\_time_{G''} \leftarrow$ The expected detection time of $(k,\tau,G'',D')$ (Eq. \ref{Exp_time})
\STATE $\Delta \leftarrow exp\_time_{G''} - exp\_time_{G'}$
\STATE $x \leftarrow$ A random value between $0$ and $1$
\IF{ $\Delta < 0$ \textbf{OR} $x<e^{\frac{-\Delta}{Temp}}$}
\STATE $G' \leftarrow G''$
\STATE $exp\_time_{G'} \leftarrow exp\_time_{G''}$
\ENDIF
\STATE $Temp \leftarrow 0.9 \times Temp$
\ENDFOR
\ENDFOR
\STATE \textbf{return} $(best\_stratey,best\_exp\_time)$
\end{algorithmic}
\label{HeuristicSearch}
\end{algorithm}

\begin{table}[]
\caption{Results of Model 1 for 48 different combinations of the input parameters. Columns $C_5$, $C_6$ and $C_{12}$ are the parameters $Max_{\tau}$, $Max_g$ and the objective value, the expected detection time, respectively. The runs with no feasible solution are shown by blank cells. }
\label{Table1}
\tiny
\begin{tabular}{|l|l|l|l|l|l|l|l|l|l|l|l|}
\hline
       & \textbf{m}  & \textbf{n}  & \textbf{p}    & \textbf{$C_5$} & \textbf{$C_6$} & \textbf{$\kappa$} & \textbf{$k$} & \textbf{$\tau$} & \textbf{$G$}                       & \textbf{$D$}                 & \textbf{$C_{12}$} \\ \hline

Run 1  & 50 & 10 & 0.05 & 4        & 30     & 9        &   &     &                       &                 &           \\ \hline
Run 2  & 50 & 10 & 0.05 & 4        & 30     & 17       &   &     &                       &                 &           \\ \hline
Run 3  & 50 & 10 & 0.05 & 4        & 22     & 9        &   &     &                       &                 &           \\ \hline
Run 4  & 50 & 10 & 0.05 & 4        & 22     & 17       &   &     &                       &                 &           \\ \hline
Run 5  & 50 & 10 & 0.05 & 7        & 30     & 9        & 2 & 5   & \{28,22\}             & \{2,5\}         & 1.7365    \\ \hline
Run 6  & 50 & 10 & 0.05 & 7        & 30     & 17       & 3 & 6   & \{16,17,17\}          & \{2,4,6\}       & 1.4355    \\ \hline
Run 7  & 50 & 10 & 0.05 & 7        & 22     & 9        & 3 & 6   & \{16,17,17\}          & \{2,4,6\}       & 1.7587    \\ \hline
Run 8  & 50 & 10 & 0.05 & 7        & 22     & 17       & 3 & 6   & \{16,17,17\}          & \{2,4,6\}       & 1.4355    \\ \hline
Run 9  & 50 & 10 & 0.1  & 4        & 30     & 9        & 3 & 3   & \{16,17,17\}          & \{1,2,3\}       & 1.0427    \\ \hline
Run 10 & 50 & 10 & 0.1  & 4        & 30     & 17       & 3 & 3   & \{16,17,17\}          & \{1,2,3\}       & 0.8692    \\ \hline
Run 11 & 50 & 10 & 0.1  & 4        & 22     & 9        & 3 & 3   & \{16,17,17\}          & \{1,2,3\}       & 1.0427    \\ \hline
Run 12 & 50 & 10 & 0.1  & 4        & 22     & 17       & 3 & 3   & \{16,17,17\}          & \{1,2,3\}       & 0.8692    \\ \hline
Run 13 & 50 & 10 & 0.1  & 7        & 30     & 9        & 3 & 3   & \{16,17,17\}          & \{1,2,3\}       & 1.0427    \\ \hline
Run 14 & 50 & 10 & 0.1  & 7        & 30     & 17       & 3 & 3   & \{16,17,17\}          & \{1,2,3\}       & 0.8692    \\ \hline
Run 15 & 50 & 10 & 0.1  & 7        & 22     & 9        & 3 & 3   & \{16,17,17\}          & \{1,2,3\}       & 1.0427    \\ \hline
Run 16 & 50 & 10 & 0.1  & 7        & 22     & 17       & 3 & 3   & \{16,17,17\}          & \{1,2,3\}       & 0.8692    \\ \hline
Run 17 & 50 & 10 & 0.2  & 4        & 30     & 9        & 2 & 2   & \{25,25\}             & \{1,2\}         & 0.8082    \\ \hline
Run 18 & 50 & 10 & 0.2  & 4        & 30     & 17       & 2 & 2   & \{25,25\}             & \{1,2\}         & 0.7005    \\ \hline
Run 19 & 50 & 10 & 0.2  & 4        & 22     & 9        & 3 & 3   & \{16,17,17\}          & \{1,2,3\}       & 1.0427    \\ \hline
Run 20 & 50 & 10 & 0.2  & 4        & 22     & 17       & 3 & 3   & \{16,17,17\}          & \{1,2,3\}       & 0.8692    \\ \hline
Run 21 & 50 & 10 & 0.2  & 7        & 30     & 9        & 2 & 2   & \{25,25\}             & \{1,2\}         & 0.8082    \\ \hline
Run 22 & 50 & 10 & 0.2  & 7        & 30     & 17       & 2 & 2   & \{25,25\}             & \{1,2\}         & 0.7005    \\ \hline
Run 23 & 50 & 10 & 0.2  & 7        & 22     & 9        & 3 & 3   & \{16,17,17\}          & \{1,2,3\}       & 1.0427    \\ \hline
Run 24 & 50 & 10 & 0.2  & 7        & 22     & 17       & 3 & 3   & \{16,17,17\}          & \{1,2,3\}       & 0.8692    \\ \hline
Run 25 & 90 & 15 & 0.05 & 4        & 30     & 15       &   &     &                       &                 &           \\ \hline
Run 26 & 90 & 15 & 0.05 & 4        & 30     & 29       &   &     &                       &                 &           \\ \hline
Run 27 & 90 & 15 & 0.05 & 4        & 22     & 15       &   &     &                       &                 &           \\ \hline
Run 28 & 90 & 15 & 0.05 & 4        & 22     & 29       &   &     &                       &                 &           \\ \hline
Run 29 & 90 & 15 & 0.05 & 7        & 30     & 15       & 4 & 6   & \{25,20,25,20\}       & \{1,3,4,6\}     & 1.4332    \\ \hline
Run 30 & 90 & 15 & 0.05 & 7        & 30     & 29       & 6 & 7   & \{13,13,13,12,23,16\} & \{1,2,3,4,5,7\} & 1.1547    \\ \hline
Run 31 & 90 & 15 & 0.05 & 7        & 22     & 15       & 6 & 7   & \{21,17,13,13,13,13\} & \{1,3,4,5,6,7\} & 1.4673    \\ \hline
Run 32 & 90 & 15 & 0.05 & 7        & 22     & 29       & 6 & 7   & \{13,13,13,22,16,13\} & \{1,2,3,4,6,7\} & 1.1548    \\ \hline
Run 33 & 90 & 15 & 0.1  & 4        & 30     & 15       & 3 & 3   & \{30,30,30\}          & \{1,2,3\}       & 0.9035    \\ \hline
Run 34 & 90 & 15 & 0.1  & 4        & 30     & 29       & 3 & 3   & \{30,30,30\}          & \{1,2,3\}       & 0.7381    \\ \hline
Run 35 & 90 & 15 & 0.1  & 4        & 22     & 15       &   &     &                       &                 &           \\ \hline
Run 36 & 90 & 15 & 0.1  & 4        & 22     & 29       &   &     &                       &                 &           \\ \hline
Run 37 & 90 & 15 & 0.1  & 7        & 30     & 15       & 3 & 3   & \{30,30,30\}          & \{1,2,3\}       & 0.9035    \\ \hline
Run 38 & 90 & 15 & 0.1  & 7        & 30     & 29       & 3 & 3   & \{30,30,30\}          & \{1,2,3\}       & 0.7381    \\ \hline
Run 39 & 90 & 15 & 0.1  & 7        & 22     & 15       & 5 & 5   & \{18,18,18,18,18\}    & \{1,2,3,4,5\}   & 1.1904    \\ \hline
Run 40 & 90 & 15 & 0.1  & 7        & 22     & 29       & 5 & 5   & \{18,18,18,18,18\}    & \{1,2,3,4,5\}   & 0.9391    \\ \hline
Run 41 & 90 & 15 & 0.2  & 4        & 30     & 15       & 3 & 3   & \{30,30,30\}          & \{1,2,3\}       & 0.9035    \\ \hline
Run 42 & 90 & 15 & 0.2  & 4        & 30     & 29       & 3 & 3   & \{30,30,30\}          & \{1,2,3\}       & 0.7381    \\ \hline
Run 43 & 90 & 15 & 0.2  & 4        & 22     & 15       &   &     &                       &                 &           \\ \hline
Run 44 & 90 & 15 & 0.2  & 4        & 22     & 29       &   &     &                       &                 &           \\ \hline
Run 45 & 90 & 15 & 0.2  & 7        & 30     & 15       & 3 & 3   & \{30,30,30\}          & \{1,2,3\}       & 0.9035    \\ \hline
Run 46 & 90 & 15 & 0.2  & 7        & 30     & 29       & 3 & 3   & \{30,30,30\}          & \{1,2,3\}       & 0.7381    \\ \hline
Run 47 & 90 & 15 & 0.2  & 7        & 22     & 15       & 5 & 5   & \{18,18,18,18,18\}    & \{1,2,3,4,5\}   & 1.1904    \\ \hline
Run 48 & 90 & 15 & 0.2  & 7        & 22     & 29       & 5 & 5   & \{18,18,18,18,18\}    & \{1,2,3,4,5\}   & 0.9391    \\ \hline
\end{tabular}
\end{table}


\begin{table}[]
\caption{Results of Model 2 for 48 different combinations of the input parameters. Columns $C_5$, $C_6$ and $C_{12}$ are the parameters $Max_{\tau}$, $Max_g$ and the objective value, the obtained minimum portion of the time for the testing process of the residents which should be allocated by the staffs, respectively. The runs with no feasible solution are shown by blank cells. }
\label{Table2}
\tiny
\begin{tabular}{|l|l|l|l|l|l|l|l|l|l|l|l|}
\hline
       & \textbf{$m$}  & \textbf{$n$}  & \textbf{$\alpha$}    & \textbf{$C_5$} & \textbf{$C_6$} & \textbf{$\kappa$} & \textbf{$k$} & \textbf{$\tau$} & \textbf{$G$} & \textbf{$D$} & \textbf{$C_{12}$} \\ \hline
Run 1  & 50 & 10 & 0.3   & 4        & 30     & 9        & 3 & 3   & \{16,17,17\}       & \{1,2,3\}     & 8.96    \\ \hline
Run 2  & 50 & 10 & 0.3   & 4        & 30     & 17       &   &     &                    &               &         \\ \hline
Run 3  & 50 & 10 & 0.3   & 4        & 22     & 9        & 3 & 3   & \{16,17,17\}       & \{1,2,3\}     & 8.96    \\ \hline
Run 4  & 50 & 10 & 0.3   & 4        & 22     & 17       &   &     &                    &               &         \\ \hline
Run 5  & 50 & 10 & 0.3   & 7        & 30     & 9        & 3 & 3   & \{16,17,17\}       & \{1,2,3\}     & 8.96    \\ \hline
Run 6  & 50 & 10 & 0.3   & 7        & 30     & 17       &   &     &                    &               &         \\ \hline
Run 7  & 50 & 10 & 0.3   & 7        & 22     & 9        & 3 & 3   & \{16,17,17\}       & \{1,2,3\}     & 8.96    \\ \hline
Run 8  & 50 & 10 & 0.3   & 7        & 22     & 17       &   &     &                    &               &         \\ \hline
Run 9  & 50 & 10 & 0.5   & 4        & 30     & 9        & 2 & 4   & \{25,25\}          & \{2,4\}       & 5.78    \\ \hline
Run 10 & 50 & 10 & 0.5   & 4        & 30     & 17       & 4 & 4   & \{12,13,12,13\}    & \{1,2,3,4\}   & 7.66    \\ \hline
Run 11 & 50 & 10 & 0.5   & 4        & 22     & 9        & 3 & 4   & \{16,17,17\}       & \{1,2,4\}     & 6.72    \\ \hline
Run 12 & 50 & 10 & 0.5   & 4        & 22     & 17       & 4 & 4   & \{12,13,12,13\}    & \{1,2,3,4\}   & 7.66    \\ \hline
Run 13 & 50 & 10 & 0.5   & 7        & 30     & 9        & 2 & 4   & \{25,25\}          & \{2,4\}       & 5.78    \\ \hline
Run 14 & 50 & 10 & 0.5   & 7        & 30     & 17       & 4 & 4   & \{12,13,12,13\}    & \{1,2,3,4\}   & 7.66    \\ \hline
Run 15 & 50 & 10 & 0.5   & 7        & 22     & 9        & 4 & 5   & \{12,13,12,13\}    & \{1,2,3,5\}   & 6.13    \\ \hline
Run 16 & 50 & 10 & 0.5   & 7        & 22     & 17       & 4 & 4   & \{12,13,12,13\}    & \{1,2,3,4\}   & 7.66    \\ \hline
Run 17 & 50 & 10 & 0.75  & 4        & 30     & 9        & 2 & 4   & \{25,25\}          & \{2,4\}       & 5.78    \\ \hline
Run 18 & 50 & 10 & 0.75  & 4        & 30     & 17       & 3 & 4   & \{16,17,17\}       & \{1,2,4\}     & 6.72    \\ \hline
Run 19 & 50 & 10 & 0.75  & 4        & 22     & 9        & 3 & 4   & \{16,17,17\}       & \{1,2,4\}     & 6.72    \\ \hline
Run 20 & 50 & 10 & 0.75  & 4        & 22     & 17       & 3 & 4   & \{16,17,17\}       & \{1,2,4\}     & 6.72    \\ \hline
Run 21 & 50 & 10 & 0.75  & 7        & 30     & 9        & 3 & 7   & \{16,17,17\}       & \{2,4,7\}     & 3.84    \\ \hline
Run 22 & 50 & 10 & 0.75  & 7        & 30     & 17       & 3 & 4   & \{16,17,17\}       & \{1,2,4\}     & 6.72    \\ \hline
Run 23 & 50 & 10 & 0.75  & 7        & 22     & 9        & 3 & 7   & \{16,17,17\}       & \{2,4,7\}     & 3.84    \\ \hline
Run 24 & 50 & 10 & 0.75  & 7        & 22     & 17       & 3 & 4   & \{16,17,17\}       & \{1,2,4\}     & 6.72    \\ \hline
Run 25 & 90 & 15 & 0.3   & 4        & 30     & 15       &   &     &                    &               &         \\ \hline
Run 26 & 90 & 15 & 0.3   & 4        & 30     & 29       &   &     &                    &               &         \\ \hline
Run 27 & 90 & 15 & 0.3   & 4        & 22     & 15       &   &     &                    &               &         \\ \hline
Run 28 & 90 & 15 & 0.3   & 4        & 22     & 29       &   &     &                    &               &         \\ \hline
Run 29 & 90 & 15 & 0.3   & 7        & 30     & 15       &   &     &                    &               &         \\ \hline
Run 30 & 90 & 15 & 0.3   & 7        & 30     & 29       &   &     &                    &               &         \\ \hline
Run 31 & 90 & 15 & 0.3   & 7        & 22     & 15       &   &     &                    &               &         \\ \hline
Run 32 & 90 & 15 & 0.3   & 7        & 22     & 29       &   &     &                    &               &         \\ \hline
Run 33 & 90 & 15 & 0.5   & 4        & 30     & 15       & 4 & 4   & \{22,23,22,23\}    & \{1,2,3,4\}   & 7.19    \\ \hline
Run 34 & 90 & 15 & 0.5   & 4        & 30     & 29       &   &     &                    &               &         \\ \hline
Run 35 & 90 & 15 & 0.5   & 4        & 22     & 15       &   &     &                    &               &         \\ \hline
Run 36 & 90 & 15 & 0.5   & 4        & 22     & 29       &   &     &                    &               &         \\ \hline
Run 37 & 90 & 15 & 0.5   & 7        & 30     & 15       & 4 & 4   & \{22,23,22,23\}    & \{1,2,3,4\}   & 7.19    \\ \hline
Run 38 & 90 & 15 & 0.5   & 7        & 30     & 29       &   &     &                    &               &         \\ \hline
Run 39 & 90 & 15 & 0.5   & 7        & 22     & 15       &   &     &                    &               &         \\ \hline
Run 40 & 90 & 15 & 0.5   & 7        & 22     & 29       &   &     &                    &               &         \\ \hline
Run 41 & 90 & 15 & 0.75  & 4        & 30     & 15       & 3 & 4   & \{30,30,30\}       & \{1,2,4\}     & 6.56    \\ \hline
Run 42 & 90 & 15 & 0.75  & 4        & 30     & 29       & 3 & 4   & \{30,30,30\}       & \{1,2,4\}     & 6.56    \\ \hline
Run 43 & 90 & 15 & 0.75  & 4        & 22     & 15       &   &     &                    &               &         \\ \hline
Run 44 & 90 & 15 & 0.75  & 4        & 22     & 29       &   &     &                    &               &         \\ \hline
Run 45 & 90 & 15 & 0.75  & 7        & 30     & 15       & 5 & 6   & \{18,18,18,18,18\} & \{1,2,3,4,6\} & 5.21    \\ \hline
Run 46 & 90 & 15 & 0.75  & 7        & 30     & 29       & 5 & 5   & \{18,18,18,18,18\} & \{1,2,3,4,5\} & 6.25    \\ \hline
Run 47 & 90 & 15 & 0.75  & 7        & 22     & 15       & 5 & 6   & \{18,18,18,18,18\} & \{1,2,3,4,6\} & 5.21    \\ \hline
Run 48 & 90 & 15 & 0.75  & 7        & 22     & 29       & 5 & 5   & \{18,18,18,18,18\} & \{1,2,3,4,5\} & 6.25    \\ \hline
\end{tabular}
\end{table}

\end{document}